\documentclass[11pt]{amsart}%book,report,article,amsart
\usepackage{amsmath,amsxtra,amssymb}
\usepackage[all]{xy}

\newcommand{\hz}{\vspace{0.5cm}}

\renewcommand{\qed}{\hspace*{\fill}$\Box$\hz\pagebreak[1]}

\newcommand{\F}{{\mathcal F}}
\newcommand{\E}{{\mathcal E}}

\newcommand{\B}{{\mathcal B}}

\newcommand{\M}{{\mathcal M}}

\newcommand{\N}{{\mathcal N}}

\newcommand{\U}{{\mathcal U}}

\newtheorem{lemma}{Lemma}[section]
\newtheorem{prop}[lemma]{Proposition}
\newtheorem{theorem}[lemma]{Theorem}
\newtheorem{cor}[lemma]{Corollary}

\newtheorem{rem}[lemma]{Remark}

\newcommand{\re}{\begin{rem}\rm}
  \newcommand{\mar}{\end{rem}}
\newtheorem{exam}[lemma]{Example}

\newtheorem{defi}[lemma]{Definition}
\begin{document}
\title{\bf On
the operator space  $UMD$  property for noncommutative \boldmath
$L_p$ \unboldmath -spaces}
%\begin{document}
%\maketitle
\author{Magdalena Musat}
\address{Department of Mathematics, 0112\\
University of California, San Diego\\
La Jolla, CA 92093-0112}
\email{mmusat@math.ucsd.edu}
\date{}
\begin{abstract}
We study the operator space $UMD$ property, introduced by Pisier
%Motivated by the classical notion of a UMD space (i.e., {\em unconditional for martingale differences}), Pisier introduced the %operator space UMD$_p$ property
in the context of noncommutative vector-valued $L_p$-spaces.
%associated to a hyperfinite (and finite) von Neumann algebra.
%We discuss basic stability properties of UMD$_p$ operator spaces.%
It is unknown whether the property is independent of $p$ in this
setting.
%We provide the first non-trivial example of a UMD$_p$\,,
%independent of $p$\,, operator space.
We prove that for $1< p, q< \infty$\,, the Schatten $q$-classes
$S_q$ are $OUMD_p$\,. The proof relies on properties of the
Haagerup tensor product and complex interpolation. Using
ultraproduct techniques, we extend this result to a large class of
noncommutative $L_q$-spaces.
Namely, we show that if $\M$ is a $QWEP$ von Neumann algebra (i.e., a quotient of a $C^*$-algebra with Lance's weak expectation property) equipped with a normal, faithful tracial state $\tau$\,, then $L_q(\M, \tau)$ is $OUMD_p$ for $1< p, q< \infty\,.$ %We discuss basic %%properties of UMD$_p$ operator spaces. We show that the UMD$_p$ property is preserved by subspaces, quotients and %ultraproducts. If an operator space $V$ is UMD$_p$\,, then its dual $V^*$ is UMD$_{p'}$\,, where $p'$ is the conjugate exponent %of $p$\,. The class of UMD$_p$ operator spaces is stable under comples interpolation.
%We provide further examples of UMD$_p$ (independent of $p$)
%operator spaces, including the non-commutative Lorentz spaces
%associated to a hyperfinite (and finite) von Neumann algebra.

\end{abstract}

\maketitle
%\input paperbmo
%\chapter{Operator Space UMD Property}
%\include{papersp}
\section{Introduction}
Probabilistic techniques are well-established powerful tools in
the study of Fourier analysis of vector-valued functions. In
particular, Banach spaces having the $UMD$ property, that is, the
property of {\em unconditionality for martingale differences} play
an important role. Deep connections with the boundedness of
certain singular integral operators, such as the Hilbert
transform, were established through the work of Burkholder,
McConnell and Bourgain. Namely, Burkholder and McConnell
\cite{Bu5} proved that if a Banach space $B$ is $UMD$\,, then the
Hilbert transform is a bounded operator on the vector-valued
Lebesgue space $L_p([0, 1]; B)$\,, for $1< p< \infty$\,. Later,
Bourgain \cite{Bo1} showed that, conversely, the boundedness of
the Hilbert transform on $L_p([0, 1]; B)$ ($1< p< \infty$) implies
that $B$ is $UMD$\,. Recall that the Banach space $B$ is $UMD$ if,
for $1< p< \infty$\,, there exists a constant $\beta_p> 0$ such
that
\begin{eqnarray*}
\left\|\sum\limits_{n=1}^k \varepsilon_ndx_n\right\|_{L_p([0, 1];
B)}&\leq& \beta_p\;\left\|\sum\limits_{n=1}^k
dx_n\right\|_{L_p([0, 1]; B)}\,,
%\label{plusminus1b}
\end{eqnarray*}
for all positive integers $k$\,, all sequences
$\varepsilon=(\varepsilon_n)_{n=1}^k$ of numbers in $\{-1, 1\}$
and all $B$-valued martingale difference sequences
$dx=(dx_n)_{n=1}^k$\,. Equivalently, for all sequences
$\varepsilon$ as above, the $\pm 1$ martingale transform
$T_\varepsilon$ generated by $\varepsilon$, i.e.,
$\,T_\varepsilon\left(\sum\limits_{n=1}^k
dx_n\right)=\sum\limits_{n=1}^k \varepsilon_ndx_n\,,$ is a bounded
operator on ${L_p([0, 1]; B)}$\,, with norm estimate
\[ \|T_\varepsilon:{L_p([0,
1]; B)}\rightarrow {L_p([0, 1]; B)}\|\leq \beta_p\,. \] The fact
that the finiteness of $\beta_p$ for some $1< p< \infty$ implies
its finiteness for all such $p$ was first proved by Pisier; see
Maurey \cite{Ma}. Results of Burkholder \cite{Bu4} provided the
first example of a $UMD$ Banach space, namely, the real line
$\mathbb{R}$\,.
%If $B$ is a UMD Banach space, the same holds for subspaces and
%quotients of $B$, for its dual $B^*$, as well as for the
%vector-valued spaces $l_p^B$ and, respectively, $L_p^B[0, 1]$, for
%all $1< p< \infty\,.$
%Aldous \cite{Al} showed that if $\,1< p<
%\infty$ and $L_p^B[0, 1]$ has an unconditional basis, then $B$ is
%UMD. Conversely, if $B$ is a UMD Banach space with an
%unconditional basis, then, for $\,1< p< \infty$\,, the space
%$L_p^B[0, 1]$ has an unconditional basis.
%The UMD property is preserved by complex, as well as real
%interpolation.
%i.e., if $(B_0, B_1)$ is a compatible couple of UMD
%Banach spaces, then the interpolation spaces $[B_0, B_1]_\theta$
%and $[B_0, B_1]_{\theta, p}$ are UMD, for $0< \theta< 1$ and $1<
%p< \infty\,.$
Other examples
%of UMD Banach spaces
include the Schatten $p$-classes for $1< p< \infty$ (Guti\'{e}rrez
\cite{Gu}, Bourgain \cite{Bo2}), and the noncommutative $L_p(\M,
\tau)$-spaces associated with a von Neumann algebra $\M$\,,
equipped with a normal, semifinite, faithful (abbreviated as
n.s.f.) trace $\tau$ (Berkson, Gillespie and Muhly \cite{BGM}). We
refer to Burkholder \cite{Bu9, Bu} for more properties of $UMD$
Banach spaces, connections to other topics and further references.
%The spaces $l_1\,, l_\infty\,, L_1[0, 1]$ and $L_\infty[0, 1]$ are not UMD.

More recently, Pisier \cite{P} developed a theory of
noncommutative vector-valued Lebesgue spaces $L_p(\M; E)$
associated with a von Neumann algebra $\M$ with an n.s.f. trace
$\tau\,.$ Two restrictions are required for the theory to be
satisfactory: $\M$ has to be hyperfinite, and $E$ equipped with an
operator space structure, that is, a sequence of matrix norms
$\|\cdot\|_m$ defined on $M_m(E)$ for each positive integer $m$\,,
such that for all $x\in M_m(E)\,, y\in M_n(E)$ and $\alpha\,,
\beta\in M_m(\mathbb{C})\,,$
\begin{eqnarray}
\|x\oplus y\|_{m+n}&=&\max\{\|x\|_m, \|y\|_n\}\,, \quad
\|{\alpha}x{\beta}\|_m\,\leq \,\|\alpha\|
\|x\|_m\|\beta\|\,.\label{ruanaxioms}
\end{eqnarray}
We recall that the class of hyperfinite von Neumann algebras
includes the algebra of bounded linear operators on a separable
Hilbert space (in particular, matrix algebras), the classical
$L_\infty$-spaces and group von Neumann algebras associated to
amenable groups, and it is closed under von Neumann algebra tensor
products. All the stability properties of the noncommutative
Lebesgue spaces $L_p(\M; E)$(e.g., duality) should be formulated
in the category of operator spaces.

Noncommutative conditional expectations and martingales arise
naturally in this setting. The $L_p$-theory of noncommutative
martingales has achieved a rapid and considerable progress in
recent years, see, e.g., Junge \cite{Ju}, Junge and Xu \cite{JX,
JX2, JX1} and Randrianantoanina \cite{Ra, Ra3}. Also,
noncommutative $BMO$ spaces were studied in \cite{Mei, Mu1, JM}.
The systematic investigation of various noncommutative martingale
inequalities started from the seminal paper \cite{PX} of Pisier
and Xu, where they introduced noncommutative Hardy spaces of
martingales and proved the analogue of the Burkholder-Gundy square
function inequalities. As a consequence, it follows that for $1<
p< \infty$\,, the $\pm 1$ martingale transform $T_\varepsilon$ is
a bounded operator on the noncommutative Lebesgue space $L_p(\M,
\tau)$\,, associated with a von Neumann algebra $\M$ with an
n.s.f. tracial state $\tau$\,. This led naturally to formulating
an appropriate notion of operator space $OUMD_p$ property in this
setting, as introduced by Pisier in \cite{P}, and to obtaining the
first example of an operator space that is $OUMD_p$ for $1< p<
\infty$\,, namely the complex plane $\mathbb{C}$\,.

%The construction of the noncommutative Lebesgue spaces, as well as
%some of their stability properties (e.g. duality) are briefly
%discussed in the Preliminaries. Therein we also explain how
%noncommutative conditional expectations and martingales arise
%naturally in this setting.

In this paper we study basic stability properties of $OUMD_p$
operator spaces and consider some related questions formulated in
\cite{P}. Our paper is organized as follows. The construction of
the noncommutative vector-valued Lebesgue spaces, as well as some
of their stability properties (e.g. duality) are briefly discussed
in the Preliminaries. Section 3 is devoted to the study of basic
stability properties of $OUMD_p$ operator spaces. Namely, we show
that, as in the classical setting, the $OUMD_p$ property is
inherited by subspaces and quotients, and it is preserved under
complex interpolation and by ultraproducts. If an operator space
$E$ is $OUMD_p$\,, then its (standard) dual $E^*$ is
$OUMD_{p'}$\,, where $p'$ is the conjugate exponent of $p$\,.
Also, each matrix level $M_m(E)$ is $OUMD_p$\,. We end Section 3
with an example (based on a construction of Pisier from \cite{P})
of a Hilbert space, subspace of some commutative $C^*$-algebra,
which is $UMD$ as a Banach space, but not $OUMD_p$\,, for any $1<
p< \infty\,.$ Section 4 contains our main results. It is unknown
whether the property is independent of $p$ in this setting. We
provide the first non-trivial example of an operator space that is
$OUMD_p$, independent of $p$\,, namely, we prove the following:
\begin{theorem}
Let $1< p, q< \infty$. Then the Schatten $q$-class $S_q$  is
$OUMD_p$\,.
\end{theorem}
The proof relies on properties of the Haagerup tensor product and
complex interpolation, and it was inspired by a question of Z.-J.
Ruan as to whether the column Hilbert space $C$ is $OUMD_p$ for
some (all) $1< p< \infty$\,.  As an application of ultraproduct
results due to Junge \cite{Ju}, it follows that a large class of
noncommutative $L_q$-spaces are ${OUMD}_p\,$, independent of
$p\,$. Namely, if $\M$ is a $QWEP$ von Neumann algebra (i.e., a
quotient of a $C^*$-algebra with the weak expectation property of
Lance \cite{La1}, see Kirchberg \cite{Ki}), equipped with an
n.s.f. tracial state $\tau$, then $L_q(\M, \tau)\,$ is $OUMD_p$
for $1< p, q< \infty\,.$ Furthermore, we show that the class of
operator spaces which are $OUMD_p$\,, independent of $p$\,,
contains all finite dimensional operator spaces, the vector-valued
Schatten classes $S_u[S_v]$ for $1< u, v< \infty$, as well as the
noncommutative Lorentz spaces associated to a hyperfinite (and
finite) von Neumann algebra. We end Section 4 with some
intermediate results towards answering Z.-J. Ruan's question,
which remains open.

\section{Preliminaries}
We refer to Effros and Ruan \cite{ER} and Pisier \cite{Pi1} for
details on operator spaces and completely bounded maps. We shall
briefly recall some definitions. A ({\em concrete}) {\em operator
space} on a Hilbert space $H$ is a norm closed linear subspace $E$
of $\B(H)$\,. For any positive integer $m$, the natural inclusion
$M_m(E)\subseteq {M_m(\B(H))}={\B(H^m)}$ induces a norm
$\|\cdot\|_m$ on $M_m(E)$\,. Ruan \cite{Ru} gave an abstract
characterization of operator spaces in terms of their matrix
norms. Namely, an ({\em abstract}) {\em operator space} is a
vector space $E$ equipped with matrix norms $\|\cdot\|_m$ on
$M_m(E)$\,, for each positive integer $m$\,, satisfying axioms
(\ref{ruanaxioms}). The morphisms in the category of operator
spaces are {\em completely bounded maps}. Given a linear map
between two operator spaces $\phi : E_0 \rightarrow E_1$\,, define
${\phi}_m : {M_m(E_0)} \rightarrow {M_m(E_1)}$ by
%\begin{eqnarray*}
${\phi}_m([v_{ij}])= [\phi(v_{ij})]$\,,
%\end{eqnarray*}
for all $ [v_{ij}]_{i,j=1}^m\in M_m(E_0)$\,. Let
$\|\phi\|_{cb}=\sup \{\|{\phi}_m\| : m\in {\mathbb{N}} \,\}$\,.
The map $\phi$ is called {\em completely bounded} if
$\|\phi\|_{cb} <\infty\,,$ and $\phi$ is called {\em completely
isometric} if all ${\phi}_m$ are isometries. The space of all
completely bounded maps from $E_0$ to $E_1$ is denoted by
$\mathcal{C}\B(E_0, E_1)$\,. Then $\mathcal{C}\B(E_0, E_1)$ is an
operator space with matrix norms defined by
\[ M_m(\mathcal{C}\B(E_0, E_1))=\mathcal{C}\B(E_0, M_m(E_1))\,, \] for
all positive integers $m$\,. The dual of an operator space $E$ is,
again, an operator space $E^*=\mathcal{C}\B(E, \mathbb{C})$\,. If
$F$ is a closed subspace of $E\,,$ then both $F$ and $E/F$ are
operator spaces; $F$ is equipped with the induced operator space
structure from $E$\,, while on $E/F$ the matrix norms are defined
by $M_m(E/F)=M_m(E)/{M_m(F)}$ for all positive integers $m$. Let
$(E_0, E_1)$ be a compatible couple of operator spaces.
%Following Xu \cite{X}\,, we equip the space $E_0+E_1$ with the
%Banach space norm $\|x\|_{E_0+E_1}=\inf\{\max(|x_0\|_{E_0},
%\|x_1\|_{E_1}): x={x_0}+{x_1}, x_0\in {E_0}, x_1\in {E_1}\}$\,.
Recall the spaces
\begin{eqnarray*}
\!\!\!\!\!\!E_0{\oplus_p} E_1&=&\{x=(x_0, x_1)\in E_0\oplus {E_1}
: \|x\|_p={(\|x\|_{E_0}^p}+{\|x_1\|_{E_1}^p)^{1/p}}\}\,, \quad 1\leq p< \infty,\\
E_0{\oplus_\infty} E_1&=&\{x=(x_0, x_1)\in E_0\oplus {E_1} :
\|x\|_\infty=\max\{\|x\|_{E_0}, \|x_1\|_{E_1}\}\}\,.
\end{eqnarray*}
As noted in \cite{Pi2}, $E_0{\oplus_\infty} E_1$ is an operator
space with matrix norms defined by
\begin{eqnarray*}
M_m(E_0{\oplus_\infty} E_1)&=&M_m(E_0)\oplus_\infty
M_m(E_1)\,,
%\label{matrixoplusinfty}
\end{eqnarray*}
for all positive integers $m$\,, while the isometric embedding
%Also, $E_0{\oplus_1} E_1$ is an operator space with matrix norms induced by the isometric embedding
\begin{eqnarray*}
\!\!\!\!\!\!\!\!\!\!\!{M_m(E_0{\oplus_1} E_1)}\,\hookrightarrow
\,{ \mathcal{C}\B(E_0^*{\oplus_\infty} {E_1^*},
M_m)}&=&M_m(E_0^*{\oplus_\infty} {E_1^*}) \,=\,M_m((E_0{\oplus_1}
E_1)^{**})
\end{eqnarray*}
induce an operator space structure on $E_0{\oplus_1} E_1$\,. For
$1< p< \infty$\,, equip $E_0{\oplus_p} E_1$ with the operator
space structure given by the isometric identification
%Note that, for $\,1< p<\infty$\,,
$E_0{\oplus_p} E_1=l_p(\{E_0, E_1\})$\,, where $l_p(\{E_0, E_1\})$
is the $l_p$-direct sum of $E_0$ and $E_1$\,. Furthermore, note
that $E_0\cap E_1$ can be identified with the diagonal
$\Delta=\{(x, x)\in E_0\oplus E_1\}$ of $E_0{\oplus_\infty}
E_1$\,, while $E_0+E_1=(E_0{\oplus_\infty} E_1)/N$\,, where
$N=\{(x_0, x_1)\in E_0\oplus E_1 : x_0+x_1=0\}\,.$ These
identifications are then used to equip $E_0\cap E_1$ and,
respectively, $E_0+ E_1$\, with appropriate operator space matrix
norms. Following Pisier \cite{Pi2}, for $0< \theta< 1$\,, we endow
the interpolation space $[E_0,E_1]_{\theta}$ with a canonical
operator space structure by defining for all positive integers
$m\,,$
\begin{eqnarray}
M_m([E_0, E_1]_{\theta})= [M_m(E_0),
M_m(E_1)]_{\theta}\,.\label{mninte}
\end{eqnarray}
Recall that the complex method of interpolation is an exact
functor of exponent $\theta$\,. Thus, if $(E_0, E_1)$ and $(F_0,
F_1)$ are two compatible couples of operator spaces, and if a map
$u: {E_0+E_1}\rightarrow {F_0+F_1}$ is completely bounded both
from $E_0$ to $E_1$ and from $F_0$ to $F_1$\,, then $u$ is
completely bounded from $[E_0, E_1]_\theta$ to $[F_0,
F_1]_\theta$\,, and, moreover, the following norm estimate holds:
\begin{eqnarray}
\|u:[E_0, E_1]_\theta\rightarrow [F_0, F_1]_\theta\|_{cb}&\leq
&\|u:E_0\rightarrow E_1\|_{cb}^{1-\theta} \|u:F_0\rightarrow
F_1\|_{cb}^\theta\,.\label{rel777}
\end{eqnarray}
%We now briefly recall the basic definitions and discuss some important aspects of the tensor product theory between operator spaces.\\[0.3cm]
Given operator spaces $E\hookrightarrow \B(H)$ and
$F\hookrightarrow \B(K)$\,, the embedding $E\otimes
F\hookrightarrow \B(H\otimes_2 K)$ induces an operator space
matrix norm $\|\cdot\|_\vee$ on $E\otimes F\,.$ It is proved in
\cite{BPa} that this matrix norm is independent on the choice of
the Hilbert spaces $H$ and $K\,.$ The completion of $E\otimes F$
with respect to the norm $\|\cdot\|_\vee$ is called the  {\em
injective tensor product} of $E$ and $F$\,. The {\em projective
tensor product} of $E$ and $F$  is defined such that the complete
isometry \[ (E\hat{\otimes} F)^*\cong{\mathcal C}\B(E, F^*)\, \]
holds. Furthermore, the {\em Haagerup tensor product} of $E$ and
$F$ is defined as the completion of $E\otimes F$ with respect to
the matrix norms
\begin{eqnarray*}
\|u\|_{h, m}&=&\inf\{\|v\|\|w\|: u=v\odot w, v\in M_{m, r}(E)\,,
w\in M_{r, m}(F)\,, r\in \mathbb{N}\}\,,
%\label{haagprnorm}
\end{eqnarray*}
where the element $v\odot w\in M_m(E\otimes F)$ is defined by
$(v\odot w)_{ij}=\sum\limits_{k=1}^m v_{ik}\otimes w_{kj}$\,, for
all $1\leq i, j\leq m\,.$ The Haagerup tensor product is both
injective and projective, associative, self-dual in the finite
dimensional case (see \cite{ER1}) and, in general, not
commutative. Moreover, it behaves very nicely with respect to
interpolation (see \cite{Ko}, \cite{Pi2}), namely,
\begin{theorem}\textup{(Kouba)}\label{koubast}
Let $(E_0, E_1)$ and $(F_0, F_1)$ be two compatible couples of
operator spaces. Then $(E_0\,{{\otimes}^{{h}}} F_0,
E_1\,{{\otimes}^{{h}}} F_1)$ is a compatible couple of operator
spaces, and for all $0< \theta< 1$ we have a complete isometry
\begin{eqnarray*}
{[E_0\,{{\otimes}^{{h}}} F_0, E_1\,{{\otimes}^{{h}}}
F_1]}_\theta&=& [E_0, E_1]_\theta{{\otimes}^{{h}}} [F_0,
F_1]_\theta\,.\label{koub}
\end{eqnarray*}
%where $E_\theta=[E_0, E_1]_\theta\,,$ respectively,
%$F_\theta=[F_0, F_1]_\theta\,.$
\end{theorem}
We refer to the Appendix in \cite{Mu} for a detailed proof, based
on ideas of Pisier from \cite{Pi3}.

The Schatten $p$-classes $S_p$ ($1\leq p\leq \infty$) are
non-commutative analogues of the Banach spaces $l_p$\,. We briefly
recall the definition and discuss their operator space structure.
If $m$ is a positive integer, denote by $S_\infty^m$ the space
$M_m\,,$ equipped with the norm $\|\cdot\|_\infty$ determined by
its identification with $\B(l_2^m)\,.$ Also, we denote by $S_1^m$
the space $\{\alpha\in M_m:
\|\alpha\|_1=\text{tr}((\alpha^*\alpha)^{1/2})< \infty\}\,.$ If
$1< p< \infty$\,, let $S_p^m=\{\alpha\in M_m: {|\alpha|}^p\in
S_1^m\}\,.$ It follows that, isometrically, $S_p^m=[S_\infty^m,
S_1^m]_{\frac1p}$\,.
%\begin{eqnarray}
%S_p^m&=&[S_\infty^m, S_1^m]_{\frac1p}\,.\label{rel345}
%\end{eqnarray}
Note that the duality $M_m^*=S_1^m$ is given by the following {\em
parallel} duality bracket
\begin{eqnarray}
\langle \,[\beta_{ij}]\,, [\alpha_{ij}]\,\rangle &=& \sum_{i,
j=1}^m \beta_{ij}\alpha_{ij}  \,\,=\,\,
\text{tr}(\beta\alpha^{\text{t}})\,,\label{paralleldual}
\end{eqnarray}
where $\alpha^{\text{t}}$ denotes the usual transposed of the
matrix $\alpha$\,.
%i.e.,
%$(\alpha^{\text{t}})_{i, j}\,=\,\alpha_{j, i}\,, 1\leq i, j\leq
%m\,.$
Thus $S_1^m=M_m^*$ has the operator space structure of the
standard dual of $M_m\,.$ Following Pisier \cite{P}, we equip
$S_p^m$ with the operator space structure (\ref{mninte}) obtained
by interpolation. Respectively, in the infinite dimensional case,
we denote by $S_\infty$ the space ${\mathcal K}(l_2)$ of compact
operators on $l_2\,,$ equipped with the operator norm.  Then
$S_\infty$ carries a natural operator space structure. Let
$S_1=S_\infty^{\,*}\,,$ equipped with the dual operator space
structure. If $1< p< \infty\,,$ we have isometrically,
\begin{eqnarray}
S_p&=&[S_\infty, S_1]_{\frac1p}\,,
\end{eqnarray}
and we equip $S_p$ with the operator space structure
(\ref{mninte}) obtained by interpolation.\\
In the following, let $E$ be an operator space. Pisier \cite{P}
constructed by interpolation the non-commutative vector-valued
Schatten $p$-classes $S_p[E]\,, 1\leq p\leq \infty\,.$ For all
$m\geq 1$\,, define
\begin{eqnarray*}
S_\infty^m[E]\,=\,S_\infty^m\check{\otimes} E\,,&&
S_\infty[E]\,=\,S_\infty\check{\otimes} E\,,
\\ S_1^m[E]\,=\,S_1^m\hat{\otimes} E\,,
&& S_1[E]\,=\,S_1\hat{\otimes} E \,,
\end{eqnarray*}
It turns out that $(S_\infty[E], S_1[E])$ (respectively,
$(S_\infty^m[E]\,, S_1^m[E])\,$) is a compatible couple for
interpolation, and, for $1< p< \infty$ and all positive integers
$m$ we define
\begin{equation}\label{spebyinterpo}
S_p^m[E]={[\,S_\infty^m[E]\,, S_1^m[E]\,]}_{\frac1p}\,,\quad
S_p[E]={[\,S_\infty[E]\,, S_1[E]\,]}_{\frac1p}\,.
\end{equation}
We equip $S_p^m[E]$ (respectively, $S_p[E]$) with the operator
space structure (\ref{mninte}). The noncommutative vector-valued
Schatten $p$-classes can be expressed in terms of the Haagerup
tensor product. Indeed, for $1\leq p\leq \infty$ and all positive
integers $m$\,, let
\begin{eqnarray}
C_p\,\,=\,\,[C, R]_{\frac1p}\,, && C_p^m\,\,=\,\,[C^m, R^m]_{\frac1p}\,, \label{cpandcpm}\\
R_p\,\,=\,\,[R, C]_{\frac1p}\,, && R_p^m\,\,=\,\,[R^m,
C^m]_{\frac1p}\,,\label{rpandrpm}
\end{eqnarray}
where $C$ and $R$ denote, respectively, the column Hilbert space
and the row Hilbert space. We should point out that we are using a
different notation from the one in \cite{P}. Namely, the space
$C_p$ is denoted therein by $C\left(\frac1p\right)$ or $C[p]$
(respectively, $C_p^m$ is denoted by $C_m\left(\frac1p\right)$ or
$C_m[p]$). Similarly, the space $R_p$ is denoted in \cite{P} by
$R\left(\frac1p\right)$ or $R[p]$ (respectively, $R_p^m$ is
denoted by $R_m\left(\frac1p\right)$ or $R_m[p]$). Using Kouba's
interpolation result, Pisier (see \cite{P}) proved that for $1\leq
p\leq \infty$ the following relations hold, completely
isometrically,
\begin{equation}\label{spehaage}
S_p[E]=C_p\otimes^{\text{h}} E\otimes^{\text{h}} R_p\,, \quad
S_p^m[E]=C_p^m\otimes^{\text{h}} E\otimes^{\text{h}} R_p^m\,.
\end{equation}
Let ${\frac1p}+{\frac1{p'}}=1\,.$ Under the parallel duality
bracket (\ref{paralleldual}) we have the following complete
isometries
\begin{equation*}
S_p[E]^{\,*}=S_{p'}[E^*]\,, \quad S_p^m[E]^{\,*}=S_{p'}^m[E^*]\,.
\end{equation*}
%\section{Operator space structure of noncommutative $L_p$-spaces}
%We now turn to the description of the appropriate operator space matrix norms on the noncommutative $L_p$-spaces. We will use Pisier's interpolation construction (\ref{mninte}).\\
Let $\M$ be a von Neumann algebra equipped with an n.f. tracial
state $\tau$\,. For $1\leq p <\infty$\,, the noncommutative
$L_p(\M, \tau)$ space is defined as the closure of $\M$ with
respect to the norm
\[ \|x\|_p=\tau((x^*x)^{\frac{p}{2}})^{\frac1p}\,. \]
The trace $\tau$ induces a canonical contractive embedding $j:
\M\rightarrow \M_*$ (where $\M_*$ denotes the unique predual of
$\M$), given by
\begin{equation}\label{bracket}
\langle j(x), y\rangle=\tau(xy)\,.
\end{equation}
With this embedding, $(\M, \M_*)$ is a compatible couple for
interpolation, and for $1\leq p< \infty$ we have the isometry
\[ L_p(\M, \tau)=[\M, \M_*]_{\frac1p}\,. \]
We now turn to the description of the appropriate operator space
matrix norms on the noncommutative $L_p$-spaces. The space
$L_\infty(\M)=\M$ carries a natural operator space structure,
since $\M$ is a $C^*$-algebra. In order to describe the operator
space structure on $L_1(\M, \tau)$ by keeping the trace duality
pairing (\ref{bracket}), we have to consider the opposite von
Neumann algebra ${\M}^{\text{op}}$\,, as explained in \cite{JR}.
Recall that ${\M}^{\text{op}}=\M$ as a vector space, but it is
endowed instead with the reversed multiplication $x\circ y=yx\,,$
for all $x, y\in \M \,.$ The algebra $\M^{\text{op}}$ carries a
natural operator space structure, with matrix norms defined for
all positive integers $m$ by
\begin{eqnarray}
\|[x_{ij}]\|_{M_m({\M}^{\text{op}})}&=&\|[x_{j\,i}]\|_{M_m(\M)}\,.\label{normop}
\end{eqnarray}
Following Junge and Ruan \cite{JR}, we define the operator space
structure on $L_1(\M, \tau)$ by
\begin{eqnarray}
L_1(\M, \tau) &\cong&  ({\M}^{\text{op}})_*\,.\label{opspaceL1}
\end{eqnarray}
The identification  (\ref{opspaceL1}) is given by the complete
isometry
\[ x\in L_1(\M, \tau)\mapsto \tau_x\in {(\M^{\text{op}})_*}\,, \]
where $\,\tau_x(y)=\tau(xy)\,,$ for all $\,y\in {\M^{\text{op}}}
\,.$ For $1<p<\infty$ and all positive integers $m$ define
\begin{eqnarray}
M_m(L_p(\M, \tau))&=&[M_m(\M), M_m(L_1({\M},
\tau))]_{\frac1p}\,.\label{operatornormlp}
\end{eqnarray}
These matrix norms verify the Ruan axioms (\ref{ruanaxioms}),
hence they determine the natural operator space structure on
$L_p(\M, \tau)$\,. Furthermore, the following Fubini-type theorem
holds isometrically,
\begin{eqnarray}
S_p^m[L_p(\M, \tau)]&=&L_p(M_m\otimes \M, \text{tr}_m\otimes
\tau)\,,\label{opsplp}
\end{eqnarray}
for any positive integer $m$ and $1\leq p\leq \infty$\,. Indeed,
we have
\begin{eqnarray*}
{\|[x_{ij}]\|}_{S_1^m[L_1(\M,
\tau)]}&=&\sup\limits_{\|[y_{ij}]\|_{M_m({\M}^{\text{op}})}\leq
1}\left\vert\sum_{i, j=1}^m\langle x_{ij}, y_{ij} \rangle
\right\vert =\sup\limits_{\|[y_{ij}]\|_{M_m({\M}^{\text{op}})}\leq
1}\left\vert\sum_{i, j=1}^m \tau(x_{ij}y_{ij})\right\vert\\
%&=&\sup\limits_{\|[z_{ij}]\|_{M_m(\M)}\leq 1}\left\vert\sum_{i,
%j=1}^m \tau(x_{ij}z_{j\,i}) \right\vert
&=&\sup\limits_{\|[z_{ij}]\|_{M_m({\M})}\leq
1}\left\vert(\text{tr}_m\otimes \tau)([x_{ij}]\cdot [z_{ji}])
\right\vert ={\|[x_{ij}]\|}_{L_1(M_m\otimes {\M},
\text{tr}_m\otimes {\tau})}\,.
\end{eqnarray*}
This proves (\ref{opsplp}) for $p=1$\,. For $p=\infty$ the
statement is clearly true. By Corollary 1.4 in \cite{P},
interpolation yields (\ref{opsplp}) for all $1\leq p\leq
\infty\,.$

In the following assume, moreover, that ${\M}$ is {\em
hyperfinite} i.e., ${\M}$ is the $w^*$-closure of an increasing
net $(\M_n)_{n\geq 1}$ of {\em finite dimensional} von Neumann
subalgebras. A celebrated theorem of Connes \cite{Con} establishes
the connection between hyperfiniteness and injectivity. Namely,
Connes proved that a von Neumann algebra $\M$ is hyperfinite if
and only if it is injective.
%another very important concept in operator algebras,
%namely {\em injectivity}.
Recall that a $C^*$-algebra $A$ is called {\em injective} if and
only if given $C^*$-algebras $B$ and $B_1$ such that $B\subseteq
B_1$ and a completely positive map $\phi: B\rightarrow A$\,, then
there exists a completely positive map $\phi_1: B_1\rightarrow A$
such that $\phi_1\!\restriction\! B=\phi\,.$ A von Neumann algebra
is called {\em injective} if it is injective as a $C^*$-algebra.
%\begin{theorem}\textup{(Connes, \cite{Con})}\label{connes}
%A von Neumann algebra $\M$ is hyperfinite if and only if it is
%injective.
%\end{theorem}
The following characterization of injective von Neumann algebras
(see Torpe \cite{Tor}) is very useful in applications. A von
Neumann algebra $\M\subseteq \B(H)$ is injective if and only if
there exists a norm $1$ projection $E: \B(H)\rightarrow \M$\,,
which is onto. The class of injective von Neumann algebras
includes $\B(H)$, and in particular, matrix algebras, the
classical $L_\infty$-spaces, group von Neumann algebras associated
to amenable groups, and it is closed under von Neumann algebra
tensor products. Also, if $\M$ is an injective von Neumann algebra
and $e$ is a projection in $\M$\,, then $e{\M}e$ is an injective
von Neumann algebra.

Let $E$ be an operator space.
%Pisier \cite{P} constructed a theory
%of vector-valued noncommutative $L_p$-spaces $L_p(\M;E)$
%associated to the hyperfinite von Neumann algebra $(\M, \tau)$ and
%the operator space $E$\,. The construction is done by
%interpolation.
Following Pisier \cite{P}, define
\begin{eqnarray}
L_1({\M}; E)&=&L_1({\M}, \tau)\hat{\otimes}E\,,\label{rel3434}
\end{eqnarray}
where the operator space structure of $\,L_1(\M, \tau)$ is given
by (\ref{opspaceL1}). Since $\M$ is a finite and hyperfinite von
Neumann algebra, it follows by results of Effros and Ruan (see
\cite{ER5} and \cite{ER6}) that the canonical inclusion
\[ v: \M\hookrightarrow (\M^{\text{op}})_*\cong L_1(\M, \tau) \]
is integral. Then, as explained in \cite{P}, the map $v\otimes
{\text{Id}_E}$ extends to a complete contraction
\begin{eqnarray}
\tilde{v}: {\M}\check{\otimes} E&\hookrightarrow &L_1({\M},
\tau)\hat{\otimes}E\,. \label{rel797}
\end{eqnarray}
Therefore $(\M\check{\otimes} E\,, L_1(\M; E))$ is a compatible
couple for interpolation. For $1< p< \infty\,,$ define
\begin{eqnarray}
L_p({\M}; E)&=&{[\,{\M}\check{\otimes} E, L_1({\M};
E)\,]}_{\frac1p}\,,\label{rel743}
\end{eqnarray}
and equip $\,L_p(\M; E)$ with the operator space matrix norms
(\ref{mninte}) obtained by interpolation. The following Fubini
theorem (see (3.6)' of \cite {P}) is a consequence of
(\ref{opsplp}). Let $1\leq p\leq \infty\,,$ then, for all von
Neumann algebras ${\N}$ with an n.f. tracial state $\phi\,,$ we
have the complete isometry
\begin{eqnarray}\label{fubth}
L_p({\M};
L_p({\N},\phi))&=&L_p({\M}\otimes{\N},\tau\otimes\phi)\,.
\end{eqnarray}
We now discuss duality in the vector-valued setting. We show that
if $\M$ is finite dimensional, then the duality results hold under
the trace duality pairing (\ref{bracket}), therefore the theory is
consistent with the scalar-valued case.
\begin{prop}\label{vectvaldualfd} Let $1\leq p< \infty$ and $p'$ be the conjugate exponent of $p$\,, i.e., ${1/p}+{1/{p'}}=1$\,.
If $\M$ is a finite dimensional von Neumann algebra equipped with
an n.f. tracial state $\tau$, then the following complete isometry
holds under the trace duality bracket (\ref{bracket})
\begin{eqnarray}
(L_{p'}(\M; E))^*&=&L_p(\M^{\textup{op}}; E^*)
\end{eqnarray}
%holds under the trace duality bracket (\ref{bracket}).
\end{prop}
\vspace*{0.3cm}\noindent{\em Proof.} Since $\M$ is finite
dimensional, it follows that $L_{p'}(\M;
E)=S_{p'}^{m_1}\oplus_{p'}\ldots \oplus_{p'} S_{p'}^{m_k}$ for
some positive integers $m_1\,, \ldots ,m_k$\,. For simplicity of
the argument, we will assume that
\begin{eqnarray*}
L_{p'}(\M; E)&=&S_{p'}^m[E]\,,
\end{eqnarray*}
where $m$ is a positive integer. We have
\begin{eqnarray*}\label{12345}
\|[x_{ij}]\|_{(L_{p'}(\M;
E))^*}&=&\sup\limits_{\|\,[y_{ij}]\,\|_{L_{p'}(\M; E)}\leq
1}\{|\langle [x_{ij}]\,, [y_{ij}] \rangle|\}=\sup\limits_{\|\,[y_{j\,i}]\,\|_{L_{p'}(\M; E)}\leq 1}\left\{\left|\sum_{i, j=1}^m x_{ij}^*(y_{j\,i})\right| \right\}\\
\nonumber &=&\|[x_{j\,i}]\|_{S_p^m[E^*]}
=\|[x_{ij}]\|_{L_p(\M^\textup{op}; E^*)}\,,
\end{eqnarray*}
where the last equality follows from the following considerations.
By (\ref{normop}), the map $\psi_m: \M\rightarrow
\M^{\textup{op}}$\,, defined by $\psi_m({[v_{ij}]}_{i,
j=1}^m)=[v_{j\,i}]_{i, j=1}^m$ is a complete isometry. Moreover,
it extends to a complete isometry $\psi_m: L_1(\M,
\tau)\rightarrow L_1(\M^{\textup{op}}, \tau)\,.$
%Then, the predual map $(\psi_m)_*:
By properties of the injective and projective tensor product, it
follows that $\psi_m\otimes \textup{Id}_{E^*}$ extends,
respectively,  to complete isometries
\[\begin{array}{rcll}
\M\check{\otimes} E^* & \stackrel{{\psi}_m\otimes \textup{Id}_{E^*}}{\longrightarrow} & \M^\textup{op}\check{\otimes} E^*\,, \\[0.3cm]
L_1(\M, \tau)\hat{\otimes} E^* & \stackrel{{\psi}_m\otimes
\textup{Id}_{E^*}}{\longrightarrow} & L_1(\M^\textup{op},
\tau)\hat{\otimes} E^*\,.
\end{array}\]
By (\ref{rel743}), interpolation with exponent $\theta=\frac1p$
shows that the map $\psi_m\otimes \textup{Id}_{E^*}: L_p(\M;
E^*)\rightarrow L_p(\M^\textup{op}; E^*)$ is a complete isometry,
which completes the argument.\qed

\section{Operator space ${OUMD}_p$ : definitions and properties}
%Before defining operator space UMD$_p$, we will sketch the background: Banach space UMD.\\[0.3cm]
Let $(\M, \tau)$ be a von Neumann algebra equipped with an n.f.
tracial state $\tau\,.$ Let $(\M_n)_{n\geq 1}$ be an increasing
filtration of von Neumann subalgebras of $\M$ such that
$\M=\left(\bigcup_n \M_n\right)^{-w^*}.$
%which generates $\M$ in the w$^*$-topology.
%The von Neumann algebra $(\M, \tau)$ and the filtration $(\M_n)_{n\geq 1}$ will remain fixed throughout.
Given a positive integer $n$\,, there is a unique normal
conditional expectation ${\E_n}:{\M} \rightarrow {\M_n}$ such that
$(\tau\!\restriction \!{\M_n})\circ {\E_n}=\E_n$ (see Takesaki
\cite{Tk}). For $1\leq p\leq \infty$\,, this extends to a norm 1
projection $\E_n: L_p(\M, \tau)\rightarrow L_p(\M_n, \tau)$
 satisfying the modular property
\[
\E_n(axb)=a{\E_n(x)}b\,,
\]
for all $a\in L_s(\M_n)\,$, $b\in L_r(\M_n)\,$ and $x\in
L_p(\M)\,$, where ${1/p}={1/r}+{1/s}\,$.

A {\em non-commutative $L_p(\M)$-martingale} relative to the
filtration $(\M_n)_{n\geq 1}$ is a sequence $x=(x_n)_{n\geq 1}$
such that $x_n \in {L_p(\M)}$ and $\E_n(x_{n+1})=x_n$\,, for all
positive integers $n$\,. We say that $x$ is a bounded
$L_p(\M)$-martingale if $\|x\|_p=\sup_n \|x_n\|_p <\infty$. The
{\em difference sequence} of $x$ is $dx=(dx_n)_{n\geq 1}$, where
$dx_n=x_n-{x_{n-1}}$, with $x_0=0.$ For $1<p<\infty$, as a
consequence of the uniform convexity of  the space $L_p(\M)$, we
can and will identify the space of all bounded
$L_p(\M)$-martingales with $L_p(\M)$ itself (see \cite{PX}, Remark
1.3).

\begin{prop}\label{tepscb}
Let $1< p< \infty\,.$ There exists $c_p>0$, depending only on
$p\,,$ such that
\[ \|T_{\varepsilon}: L_p(\M, \tau)\rightarrow L_p(\M, \tau)\|_{cb}\leq c_p \,, \]
where $\varepsilon$ denotes a sequence $(\varepsilon_n)_{n\geq 1}$
of numbers in $\{-1, 1\}$ and $T_\varepsilon$ is the $\pm 1$
martingale transform generated by $\varepsilon$, i.e.,
$\,T_\varepsilon\left(\sum\limits_{n=1}^k
dx_n\right)=\sum\limits_{n=1}^k \varepsilon_ndx_n\,,$ for all
positive integers $k$ and all $L_p(\M)$-martingale difference
sequences $dx=(dx_n)_{n=1}^k$ relative to the filtration
$(\M_n)_{n\geq 1}\,.$
\end{prop}
\vspace*{0.3cm} \noindent {\em Proof.} As a consequence of the
Pisier-Xu noncommutative version of the Burkholder-Gundy square
function inequalities, there exists a constant $c_p> 0$ such that
\begin{eqnarray}
\|T_{\varepsilon}: L_p(\M, \tau)\rightarrow L_p(\M, \tau)\|&\leq&
c_p \,.\label{tepsilon}
\end{eqnarray}
By Lemma 1.7 of \cite{P}, we have
\begin{eqnarray}
\|T_\varepsilon\|_{cb}&=&\sup_m\|\,\text{Id}_{S_p^m}\otimes
T_\varepsilon: S_p^m[L_p(\M, \tau)]\rightarrow S_p^m[L_p(\M,
\tau)]\,\|\,.
\end{eqnarray}
Let $m\geq 2$. We will prove that
\begin{eqnarray}
\|\text{Id}_{M_m}\otimes T_\varepsilon: S_p^m[L_p(\M,
\tau)]\rightarrow S_p^m[L_p(\M, \tau)]\|&\leq&
c_p\,.\label{optepsilon}
\end{eqnarray}
By Fubini's theorem (\ref{opsplp}) we have the isometry
$S_p^m[L_p(\M, \tau)]=L_p(M_m\otimes \M, \text{tr}_m\otimes
\tau)$\,, where $\text{tr}_m$ is the standard normalized trace on
$M_m\,.$ Note that $(M_m\otimes \M_n)_{n\geq 1}$ is a filtration
of the algebra $M_m(\M)=M_m\otimes \M\,.$ For all positive
integers $n$\,, denote $\text{Id}_{M_m}\otimes \E_n$ by
$\mathbb{E}_n\,.$ Then $\mathbb{E}_n:M_m(\M)\rightarrow M_m(\M_n)$
is the unique trace preserving conditional expectation onto
$\,M_m(\M_n)\,$. Moreover, for all $x=[x_{ij}]_{i, j=1}^m\in
M_m(L_p(\M, \tau))$ and all positive integers $n$ we have
\begin{eqnarray}
{[\,\E_n(x_{ij})\,]}_{i, j=1}^m&=&\mathbb{E}_n(x)\,.\label{rel121}
\end{eqnarray}
By applying (\ref{tepsilon}), together with (\ref{rel121}) to the
algebra $M_m(\M)$ and its filtration $(M_m\otimes \M_n)_{n\geq
1}\,,$ we obtain for all positive integers $k$ and all sequences
$\varepsilon=(\varepsilon_n)_{n=1}^k$ of numbers in $\{-1, 1\}$
\begin{eqnarray*}
\left\|\sum\limits_{n=1}^k (\text{Id}_{M_m}\otimes {\varepsilon_n
(\E_n-{\E_{n-1}})})(x)\right\|_p&\leq
&c_p\;\left\|\sum\limits_{n=1}^k (\text{Id}_{M_m}\otimes
(\E_n-{\E_{n-1}}))(x)\right\|_p\,.
\end{eqnarray*}
This shows that (\ref{optepsilon}) holds and the proof is
complete.\qed

\re\label{tepssa}\rm Let $1< p< \infty$ and
$\varepsilon=(\varepsilon_n)_{n\geq 1}$ be a sequence of numbers
in $\{-1, 1\}$\,. Then the $\pm 1$ martingale transform
$T_{\varepsilon}$ generated by $\varepsilon$ is a self adjoint
operator on $L_p(\M, \tau)\,,$ under the trace duality bracket
(\ref{bracket}). The proof of this fact is similar to the one in
the classical setting, due to Burkholder \cite{Bu4}. The key point
is that $\left(\,{\bigcup_n L_p(\M_n,
\tau_n)}\,\right)^{-\|\cdot\|_p}=L_p(\M, \tau)\,.$ By (\ref
{tepsilon})\,, $T_\varepsilon$ is a bounded linear operator on
$L_p(\M, \tau)$\,. Therefore, it suffices to prove that the
restriction of $T_\varepsilon$ to ${L_p(\M_n, \tau_n)}$ is self
adjoint, for all for all positive integers $n$\,. Indeed, if $x\in
L_p(\M_n, \tau_n)$\,, then $x=\sum_{k=1}^n d_k(x)\,,$ where
$d_k=\E_k-\E_{k-1}\,,$ with $\E_0=0\,.$
%Then, for $y=\sum\limits_{k=1}^n d_k(y)\in L_{p'}(\M_n, \tau_n)\,,$ we have  by the duality (\ref{rel6677})
%\[ \langle T_\varepsilon(x), y\rangle \,=\, \sum_{k=1}^n \langle \varepsilon_k d_k(x), d_k(y)\rangle
%\,=\,\sum_{k=1}^n \tau(\varepsilon_k d_k(x) d_k^*(y))
%\,=\,\sum_{k=1}^n \langle d_k(x), \varepsilon_x d_k(y)\rangle
%&=&\tau\left(\left(\sum_{k=1}^n d_k(x)\right)\left(\sum_{k=1}^n \varepsilon_k d_k^*(y)\right)\right)\\
%\,=\,\langle x, T_\varepsilon(y)\rangle\,, \]
Note that, for all positive integers $j$\,, $\E_j$ is a
self-adjoint operator on $L_p(\M, \tau)\,,$ since $\E_j$ is the
dual map of the canonical isometric embedding of $L_1(\M_n, \tau)$
into $L_1(\M, \tau)$\,. Therefore each $d_k$ is a self-adjoint
operator with respect to the trace duality pairing
(\ref{bracket}).  Since $T_\varepsilon(x)=\sum_{k=1}^n
\varepsilon_k d_k(x)$\,, the conclusion follows.\mar
%and we have
%\begin{eqnarray*}
%T_\varepsilon(x)&=&\sum_{k=1}^n \varepsilon_k d_k(x)\,,
%\end{eqnarray*}
%which yields the conclusion. \mar

For the remainder of this section we will assume, moreover, that
$\M$ is {\em hyperfinite}\,. Note that, consequently, each von
Neumann subalgebra $\M_n$ is hyperfinite, as well. Indeed, since
$\M$ is injective, there exists a norm $1$ projection $\E:
\B(H)\rightarrow \M$\,, which is onto. Composing $\E$ with the
conditional expectation ${\E_n}:{\M} \rightarrow {\M_n}$ we obtain
a norm $1$ projection $\E_n\circ \E: \B(H)\rightarrow \M_n$\,,
which is onto. This ensures that $\M_n$ is injective, or
equivalently, by Connes' theorem, $\M_n$ is hyperfinite.

\begin{prop} \label{vectvalcondexp} Let $\,1\leq p\leq \infty\,.$ Then, for all positive integers $\,n$\,, the conditional expectation
$\,\E_n: L_p(\M, \tau)\,\rightarrow \,L_p(\M_n, \tau_n)\,$ extends
to a complete contraction
\begin{eqnarray}
\E_n \otimes {\textup{Id}_{E}}:L_p(\M; E)&\rightarrow &L_p(\M_n;
E)\,.\label{rel1446}
\end{eqnarray}
\end{prop}
\vspace*{0.3cm} \noindent {\em Proof.} We first show that we have
the complete contraction
\begin{eqnarray}
\E_n\otimes {\textup{Id}_{E}}: \M\check{\otimes} E&\rightarrow
&\M_n\check{\otimes} E\label{rel1122}
\end{eqnarray}
By the injectivity property of the injective tenor product, it
suffices to prove that the map $\E_n: \M\rightarrow
 \M_n$ is a complete contraction. Since $\|\E_n\|=1$\,, this is
equivalent to showing that $\E_n$ is completely positive (see
Paulsen \cite{Pau}). This is, indeed, the case, since for all
$m\geq 1$\,, the map ${\textup{Id}}_{M_m}\otimes \E_n: {M_m\otimes
\M}\,\,\rightarrow \,\,{M_m\otimes \M_m}$ is the unique
trace-preserving conditional expectation onto $M_m\otimes \M_m$\,,
hence it is positive. We now show that we have the complete
contraction
\begin{eqnarray}
\E_n\otimes {\textup{Id}_{E}}: L_1(\M, \tau)\hat{\otimes}
E&\rightarrow &L_1(\M_n, \tau_n)\hat{\otimes} E\label{rel1123}
\end{eqnarray}
Since the projective tensor product is projective, it suffices to
prove that the map $\E_n: L_1(\M, \tau)\rightarrow L_1(\M_n,
\tau_n)$ is a complete contraction. Let us denote for the moment
this map by $u_n$\,, in order to avoid confusion. Note that, under
the trace duality bracket (\ref{bracket}), the dual map $u_n^*$ is
exactly $\E_n:\M \rightarrow M_n\,.$ Indeed, for $x\in L_1(\M,
\tau)$ and $y\in L_1(\M_n, \tau_n)$ we have
\begin{eqnarray}
\langle u_n(x)\,, y \rangle&=&\langle x\,, \E_n(y)
\rangle\,.\label{56}
\end{eqnarray}
It follows that $\|u_n\|_{cb}=\|u_n^*\|_{cb}=\|\E_n\|_{cb}\leq
1$\,, which proves the assertion.
%i.e., $\,\E_n:L_1(\M, \tau)\rightarrow L_1(\M_n, \tau_n)\,$ is a complete contraction, as wanted.\\
By (\ref{rel1122}) and (\ref{rel1123})\,, interpolation with
exponent $\theta=1/p$ yields the conclusion for $1< p<
\infty$\,.\qed

\re\rm In view of Proposition \ref{vectvalcondexp} we can consider
vector-valued noncommutative martingales in this setting. Note
that, as in the scalar-valued case, any $L_p(\M; E)$-martingale
with respect to the filtration $(\M_n)_{n\geq 1}$ can be
approximated by finite $L_p(\M; E)$-martingales (with respect to
the same filtration). Therefore, it suffices to consider finite
martingales only. \mar

\begin{defi}\rm\label{defumdp}
Let $E$ be an operator space and $1< p< \infty\,.$ We say that $E$
is $OUMD_p$ {\em with respect to the filtration $(\M_n)_{n\geq
1}$} if there exists a constant $c_p > 0$ such that
\begin{eqnarray}
\left\|\,\sum_{n=1}^k\varepsilon_n dx_n\,\right\|_{L_p({\M};
E)}&\leq& c_p\,\left\|\,\sum_{n=1}^k dx_n\, \right\|_{L_p({\M};
E)}\label{defumd}
\end{eqnarray}
for all positive integers $k\,,$ all sequences
$\varepsilon=(\varepsilon_n)_{n=1}^k$ of numbers in $\{-1, 1\}$
and all
 martingale difference sequences $dx\,=\,(d_n)_{n=1}^k\subset L_p({\M}; E)\,,$ relative to the filtration $(\M_n)_{n\geq 1}\,.$

If this holds for all hyperfinite von Neumann algebras ${\M}$,
equipped with an n.f. tracial state $\tau$\,, and all filtrations
of $\M$\,, we say that $E$ is $OUMD_p$\,.
\end{defi}

\re By Definition \ref{defumdp}, $E$ is $OUMD_p$ {\em with respect
to the filtration} $(\M_n)_{n\geq 1}$ if and only if there exists
a constant $c_p> 0\,,$ depending on $p$ and $E$, such that
\begin{eqnarray}
\|\, T_\epsilon\otimes {\text{Id}_E}: L_p(\M; E)\rightarrow
L_p(\M; E)\,\|&\leq& c_p\,,\label{tepsfiltr}
\end{eqnarray}
for all finite sequences $\varepsilon=(\varepsilon_n)_{n\geq 1}$
of numbers in $\{-1, 1\}$\,. Martingale differences are considered
with respect to the filtration $(\M_n)_{n\geq 1}\,.$ Note that a
priori the constant $c_p$ might also depend on the algebra $\M$
and on the filtration.
%The following result shows that if $E$ is $OUMD_p$\,, then the constant $C_p$ can be chosen independent on the algebra $\M$ %and on the filtration.
\mar

\begin{lemma}\label{constant}
If $E$ is $OUMD_p\,,$ there exists $c_p> 0$\,, depending only on
$p$ and $E$, such that
\begin{eqnarray}
\|\, T_\epsilon\otimes {\text{Id}_E}: L_p(\M; E)\rightarrow
L_p(\M; E)\,\|&\leq& c_p\,,\label{teps}
\end{eqnarray}
for all hyperfinite von Neumann algebras $(\M, \tau)\,,$ equipped
with an n.f. tracial state $\tau$\,, all filtrations of $\M$ and
all finite sequences $\varepsilon=(\varepsilon_n)_{n\geq 1}$ of
numbers in $\{-1, 1\}$\,.
\end{lemma}
\vspace*{0.3cm} \noindent {\em Proof.} Assume that there exists a
sequence $(\M_k, \tau_k)_{k\geq 1}$ of hyperfinite von Neumann
algebras, equipped with n.f. tracial states $\tau_k$\,, such that
the sequence of corresponding constants $c_p^{(k)}$ converges to
$\infty$\,. We will show that this leads to a contradiction. Let
$\M=\oplus_k M_k$ be the direct sum of the algebras $\M_k$\,. Note
that $\M$ is a hyperfinite von Neumann algebra and that we can
define an n.f. tracial state $\tau$ on $\M$ by \[
\tau((x_k)_{k\geq 1})=\sum\limits_{k\geq
1}\frac1{2^k}\,\tau_k(x_k)\,. \] For each positive integer $k$,
there is an increasing net of finite dimensional algebras
$(\M_k^{(n)})_{n\geq 1}$ whose union generates $\M_k$ in the
$w^*$-topology. For a finite sequence
$\varepsilon=(\varepsilon_n)_{n\geq 1}$ of numbers in $\{-1,
1\}$\,, let $T_\varepsilon^{(k)}$ denote the vector-valued $\pm 1$
martingale transform generated by $\varepsilon$, where martingale
differences are considered with respect to the algebra $\M_k$ and
its filtration $(\M_k^{(n)})_{n\geq 1}\,.$
%For simplicity of notation, in this proof  we will omit the superscript $\,E$ in the notation of the martingale transforms.
For all positive integers $n$\,, let
\[ \M^{(n)}\,=\,\oplus_k \M_k^{(n)}\,. \]
Then  $(\M^{(n)})_{n\geq 1}$ is a filtration of $\M$\,. Let
$T_\varepsilon$ denote the $\pm 1$ martingale transform generated
by $\varepsilon$, where martingales are considered with respect to
the algebra $\M$ and its filtration  $(\M^{(n)})_{n\geq 1}\,.$ By
our assumption it follows that for all positive integers $k$\,,
there exists $y_k\in L_p(\M_k; E)$ such that $\|y_k\|_{L_p(\M_k;
E)}=1$ and \[ \|\,(T_\varepsilon^{(k)}\otimes
{\text{Id}_E})(y_k)\,\|_{L_p(\M_k; E)}\geq {c_p^{(k)}}\,. \] Let
$Y_k=(0, \ldots, 0, \underbrace{y_k}_{k^{\text{th}}} \,, 0,\ldots)
\,.$ It follows that
\begin{eqnarray*}
\|\,(T_\varepsilon\otimes {\text{Id}_E})(Y_k)\,\|_{L_p(\M_k; E)}&=& \left(\frac1{2^k}\tau_k(\,|\,(T_\varepsilon^{(k)}\otimes {\text{Id}_E})(y_k)\,|^p\,)\right)^{\frac1p}\\
&=&\frac1{2^{k/p}}\,\,\|\, (T_\varepsilon^{(k)}\otimes
{\text{Id}_E})(y_k)\,\|_{L_p(\M_k; E)}\geq
\frac1{2^{k/p}}\,\,c_p^{(k)}\,.
\end{eqnarray*}
Moreover, we have
\begin{eqnarray*}
\|\,Y_k\,\|_{L_p(\M;
E)}=\left(\frac1{2^k}\,\tau_k(\,|y_k|^p\,)\right)^{\frac1p}
=\frac1{2^{k/p}}\,\,\|y_k\|_{L_p(\M_k; E)} =\frac1{2^{k/p}}\,.
\end{eqnarray*}
This shows that $\|\,T_\varepsilon\otimes {\text{Id}_E}: L_p(\M;
E)\rightarrow L_p(\M; E)\,\|\geq c_p^{(k)}$\,, for all positive
integers $\,k\,.$ This implies that $\|\,T_\varepsilon\otimes
{\text{Id}_E}: L_p(\M; E)\rightarrow L_p(\M; E)\,\|=\infty\,,$
%Hence, if the assertion made in
%Lemma \ref{constant} is not true, then there is a hyperfinite von
%Neumann algebra $(\M, \tau)$ equipped with an n.f. tracial state
%$\tau\,,$ and a filtration $(\M^{(n)})_{n\geq 1}\,$ of $\M$ for
%which the corresponding constant is $\infty$.
which contradicts the assumption that $E$ is $OUMD_p$\,. \qed

We denote by $c_p(E)$ the smallest constant $c_p> 0$ satisfying
(\ref{teps}).

\re\label{charactumdp} A similar argument as in the proof of
Proposition \ref{tepscb} shows that an operator space $E$ is
$OUMD_p$\,, for some $1< p< \infty$ if and only if
\begin{eqnarray}
\|\,T_\varepsilon\otimes{\text{Id}_E}: L_p(\M; E)\rightarrow
L_p(\M; E)\,\|_{cb}&\leq & c_p(E)\,,\label{defumdpconst}
\end{eqnarray}
for all hyperfinite von Neumann algebras $(\M, \tau)$ with an n.f.
tracial state $\tau$, all filtrations of $\M$ and all finite
sequences $\varepsilon=(\varepsilon_n)_{n\geq 1}$ of numbers in
$\{-1, 1\}$\,. \mar

\re \label{constBurk} If we assume that $\M$ is a commutative von
Neumann algebra in the Definition \ref{defumdp}, then we recover
the classical notion of a $UMD$ Banach space. In particular, if an
operator space $E$ is $OUMD_p$\,, for some $1< p< \infty$\,, then
$E$ is $UMD$ (as a Banach space) and, moreover,
\[ \beta_p(\mathbb{C})\leq \beta_p(E)\leq c_p(E)\,. \]
It was proved by Burkholder \cite{Bu1} that
$\beta_p(\mathbb{C})=p^*-1$\,, where $p^*=\max\{p, {p}/(p-1)\}$\,.
 \mar

\re\label{Complumdp} By Proposition \ref{tepscb} and Remark
\ref{charactumdp}, the operator space $E=\mathbb{C}$ is
$OUMD_p\,,$ for $1< p< \infty\,.$ Moreover, as proved by
Randrianantoanina \cite{Ra}, the corresponding constant
$c_p(\mathbb{C})$ has the same optimal order of growth as
$\beta_p(\mathbb{C})$\,.
%Namely, $c_p(\mathbb{C})\approx O(p)$ as $p\rightarrow \infty\,.$
This estimate follows directly by interpolation from a striking
result proved in \cite{Ra}, namely the fact that the
noncommutative $\pm 1$ martingale transforms are of weak type $(1,
1)\,.$ \mar

\re \label{umdpislocal} $OUMD_p$ is a local property, i.e., if
there exists an increasing sequence of closed subspaces
$(E_k)_{k\geq 1}$ of $E$\,, whose union generates $E$\,, such that
each $E_k$ is $OUMD_p$ and the corresponding constants satisfy
$\sup_k c_p(E_k)< \infty\,,$ then $E$ is $OUMD_p$\,. \mar

The following proposition summarizes some of the properties of
$OUMD_p$ operator spaces. We refer to Pisier \cite{Pi1, P} for
background on the ultraproduct theory for operator spaces.

\begin{prop}\label{propertiesumdp}
Let $1< p< \infty$ and $E$ be an operator space.
\begin{enumerate}
\item If $E$ is $OUMD_p$ and $F$ is an operator space completely
isomorphic to $E\,,$ then $F$ is $OUMD_p\,.$ \item If $E$ is
$OUMD_p$ and $F$ is a closed subspace of $E$, then both $F$ and
$E/F$ are $OUMD_p$\,. \item If $E$ is $OUMD_p\,,$ then $L_p(\N;
E)$ is $OUMD_p$\,, for all hyperfinite von Neumann algebras
$\N\,,$ equipped with an n.f. tracial state $\tau\,.$ \item If $E$
is $OUMD_p$\,, then $M_m(E)$ is $OUMD_p$\,, for all positive
integers $m$\,. \item If $E$ is $OUMD_p$\,, then its dual $E^*$ is
$OUMD_{p'}\,,$ where ${1/p}+{1/{p'}}=1\,.$ \item Let $I$ be an
index set and $(E_i)_{i\in I}$ be a family of $OUMD_p$ operator
spaces. Assume that the corresponding $OUMD_p$ constants satisfy
$\sup_i c_p(E_i)<\infty\,.$ Let $\mathcal U$ be an ultrafilter on
$I$\,. Then the ultraproduct $\hat{E}=(E_i)_{\U}$ is $OUMD_p$\,.
\item If $E$ and $F$ are $OUMD_p$\,, then $E\oplus_q F$ is
$OUMD_p$\,, for $1< q< \infty$\,. \item Let $1<q, s< \infty$ and
$0< \theta< 1$ be such that
${({1-\theta})/{q}}+{{\theta}/{s}}=1/p\,.$ If $(E_0, E_1)$ is a
compatible couple of operator spaces such that $E_0$ is $OUMD_q$
and $E_1$ is $OUMD_s\,,$ then $E_\theta=[E_0, E_1]_\theta$ is
$OUMD_p$\,.
\end{enumerate}
\end{prop}
\vspace*{0.3cm} \noindent {\em Proof.} Let $c_p(E)$ denote the
$OUMD_p$ constant of $E$\,. Throughout the proof, let $(\M, \tau)$
be a hyperfinite von Neumann algebra with an n.f. tracial state
$\tau$\,, and $(\M_n)_{n\geq 1}$ a filtration of $\M$\,.
Martingale differences will be considered with respect to this
filtration. Let $\varepsilon=(\varepsilon_n)_{n\geq 1}$ be a
finite sequence of numbers in $\{-1, 1\}$ and denote by
$T_\varepsilon$ the corresponding $\pm 1$ martingale transform.
For each statement, we will prove the boundedness of the
appropriate {\em vector-valued} martingale transform.\\
(1)
%\begin{eqnarray}
%\|\,T_\varepsilon\otimes {\text{Id}_E}: L_p(\M; E)\rightarrow
%L_p(\M; E)\,\|_{cb}&\leq & c_p(E)\,. \label{rel98}
%\end{eqnarray}
Let $u$ be a complete isomorphism $u: E\rightarrow F\,.$  Then, as
observed in (3.1) of \cite{P}, the map $\text{Id}_{L_p(\M,
\tau)}\otimes u$ extends to a complete isomorphism $\tilde{u}:
L_p(\M; E)\rightarrow L_p(\M; F)$\,. Hence, there exists a
constant $c_p^{\prime}(F)\,,$ such that $\|\,T_\varepsilon\otimes
{\text{Id}_F}: L_p(\M; F)\rightarrow L_p(\M; F)\,\|_{cb}\leq
{c'}_p(F)\,.$
By Remark \ref{charactumdp}, this implies that $F$ is $OUMD_p$\,.\\
(2) By (3.4) of \cite{P}\,, the space $L_p(\M; F)$ can be
identified with a closed subspace of $L_p(\M; E),$ and we have a
complete isometry $L_p(\M; E/F){\cong } {{L_p(\M; E)}/{L_p(\M;
F)}}\,.$ Hence
\begin{eqnarray*}
\|\,T_\varepsilon\otimes {\text{Id}_F}: L_p(\M; F)\rightarrow
L_p(\M; F)\,\|_{cb}&\leq & c_p(E)\,,\\
\|\,T_\varepsilon\otimes {\text{Id}_{E/F}}: L_p(\M;
E/F)\rightarrow L_p(\M; E/F)\,\|_{cb}&\leq & c_p(E)\,.
\end{eqnarray*}
Therefore $F$ and $E/F$ are both $OUMD_p\,.$ Moreover, the
corresponding $OUMD_p$ constants satisfy
\begin{equation}
c_p(F)\,\leq \,c_p(E)\,, \quad c_p(E/F)\,\leq
\,c_p(E)\,.\label{rel778}
\end{equation}
(3) An application of Fubini's theorem (\ref{fubth}) yields the
complete isometry
\begin{eqnarray}
L_p(\M; L_p(\N; E))&\cong&L_p(\M\otimes \N; E)\,.\label{rel39}
\end{eqnarray}
The assertion follows from the fact that $E$ is $OUMD_p$ with
respect to the hyperfinite algebra $\M\otimes \N$ and its
filtration $(\M_n\otimes \N)_{n\geq 1}$\,. Moreover, the following
estimate holds
\begin{eqnarray}
c_p(L_p(\N; E))&\leq &c_p(E)\,.\label{rel781}
\end{eqnarray}
(4) We prove that $M_m(E)=M_m\check{\otimes} E$ is completely
isomorphic to $S_p^m[E]=L_p(M_m; E)\,.$ Equivalently, we show that
$d_{cb}(\,M_m(E), L_p(M_m; E)\,)< \infty\,,$ where the {\em
c.b.-Banach-Mazur distance} between two operator spaces $E_0$ and
$E_1$ is defined by
\begin{eqnarray*}
d_{cb}(E_0, E_1)&=&\inf\{\|\phi\|_{cb}\|\phi^{-1}\|_{cb}: \phi:
E_0\rightarrow E_1 \,\,\mbox{is an isomorphism}\,\}\,.
\end{eqnarray*}
Recall, by (\ref{rel797}) the complete contraction $\tilde{v}:
M_m(E)\hookrightarrow L_1(M_m)\hat{\otimes} E$\,, induced by the
canonical inclusion $v: \M\hookrightarrow  L_1(\M, \tau)$\,. By
interpolation it follows that
\begin{eqnarray}
\|\tilde{v}: M_m(E)\rightarrow L_p(M_m; E)\|_{cb}&\leq &1\,.
\end{eqnarray}
It remains to show that ${\tilde{v}}^{-1}\in
\mathcal{C}\B(L_p(M_m; E)\,, M_m(E))\,.$ By results of Effros and
Ruan \cite{ER}, the map ${\tilde{v}}^{-1}: L_1(M_m)\hat{\otimes}
E\rightarrow M_m\check{\otimes} E$ is completely bounded if and
only if the map $v^{-1}: L_1(M_m)\rightarrow M_m$ is nuclear, and,
moreover,
\begin{eqnarray}
\|\,{\tilde{v}}^{-1}: L_1(M_m)\hat{\otimes} E\rightarrow
M_m\check{\otimes} E\,\|_{cb}&\leq& \nu(v^{-1})\,,\label{rel541}
\end{eqnarray}
 where $\nu(v^{-1})$ denotes the nuclear norm of $v^{-1}$\,. Let $\{e_{i, j}, e_{i, j}^*\}_{1\leq i, j\leq m}\in M_m\times M_m^*\,$ be an Auerbach basis for $M_m\,.$ Then
\begin{eqnarray*}
\nu(v^{-1})\leq \left\|\sum_{i, j=1}^m e_{i, j}^*(e_{i,
j})\right\| \leq \sum_{i, j=1}^m \left\|e_{i, j}^*(e_{i,
j})\right\|\,\leq \,m^2\,.
\end{eqnarray*}
From (\ref{rel541}) and the fact that ${\tilde{v}}^{-1}$ is a
complete contraction on $M_m(E)$\,, we obtain by interpolation
\begin{eqnarray}
\|\,{\tilde{v}}^{-1}: L_p(M_m; E)\rightarrow M_m(E)\,\|_{cb}&\leq
&m^{\frac{2}{p}}\,.
\end{eqnarray}
This implies that
\begin{eqnarray*}
d_{cb}(\,M_m(E), L_p(M_m; E)\,)&\leq &m^{\frac{2}{p}}\,<
\,\infty\,.
\end{eqnarray*}
The assertion follows as a consequence of Items (1) and (3)\,.\\
%By properties (1) and (3) it follows that $M_m(E)$ is UMD$_p\,.$\\
(5) Since $\M$ is hyperfinite, there exists a filtration
$(\N_\alpha)_{\alpha\geq 1}$ of {\em finite dimensional}
subalgebras whose union generates $\M$\,. For ${\alpha\geq 1}$\,,
denote by $\tilde{\E_\alpha}$ the unique trace-preserving
conditional expectation from $\M$ onto $\N_\alpha\,.$ By Remark
3.2 in \cite{P} (see also Theorem 3.4. therein) it follows that
$\left(\,\bigcup_\alpha L_p(\N_\alpha;
E^*)\,\right)^{-\|\cdot\|_p}=L_p(\M; E^*)$\,. This implies that
\begin{eqnarray}
\|\,T_\varepsilon\otimes {\text{Id}_{E^*}}\,\|_{cb}&\leq &
\sup_\alpha \|\,\tilde{\E_\alpha}(T_\varepsilon\otimes
{\text{Id}_{E^*}})\tilde{\E_\alpha}: L_{p'}(\N_\alpha; E^*)
\rightarrow L_{p'}(\N_\alpha; E^*)\,\|_{cb}\,.\label{relnicevstar}
\end{eqnarray}
The following complete isometry holds under trace duality (see
Proposition \ref{vectvaldualfd}):
\begin{eqnarray}
L_{p'}(\N_\alpha; E^*)&=&(L_p(\N_\alpha^{\,\textup{op}};
E))^{\,*}\,.\label{rel87960}
\end{eqnarray}
Furthermore, since $\left(\,\bigcup_\alpha
L_p(\N_\alpha^{\,\textup{op}};
E)\,\right)^{-\|\cdot\|_p}=L_p(\M^{\,\textup{op}}; E)$ and $E$ is
$OUMD_p$ with respect to $\M^{\,\textup{op}}$, it follows that
\begin{eqnarray*}
\|\,T_\varepsilon\otimes {\text{Id}_{E}}\,\|_{cb}&\leq &
\sup_\alpha \|\,\tilde{\E_\alpha}(T_\varepsilon\otimes
{\text{Id}_{E}})\tilde{\E_\alpha}: L_p(\N_\alpha^\textup{op}; E)
\rightarrow L_p(\N_\alpha^\textup{op}; E)\,\|_{cb}\,\,\leq
\,\,c_p(E)\,.
%\label{relnicev}
\end{eqnarray*}
By passing to the dual and using the fact that the conditional
expectations $\tilde{\E_\alpha}$ and $T_\varepsilon$ are self-dual
under trace duality (see Remark \ref{tepssa}), we obtain
\begin{eqnarray}
\|\,\tilde{\E_\alpha}(T_\varepsilon\otimes
{\text{Id}_E})^*{\tilde{\E_\alpha}}: (L_p(\N_\alpha^\textup{op};
E))^{\,*}\rightarrow (L_p(\N_\alpha^\textup{op};
E))^{\,*}\,\|_{cb}&\leq & c_p(E)\,.\label{rel40}
\end{eqnarray}
Applying (\ref{rel87960}) together with (\ref{relnicevstar}), it
follows that
\begin{eqnarray}
\|\,T_\varepsilon\otimes {\text{Id}_{E^*}}: L_{p'}(\M;
E^*)\rightarrow L_{p'}(\M; E^*)\,\|_{cb}&\leq &
c_p(E)\,.\label{rel80900}
\end{eqnarray}
Hence $E^*$ is $OUMD_{p'}$\,, and we have the following estimate
for the corresponding constant
\begin{eqnarray}
c_{p'}(E^*)&\leq &c_p(E)\,.\label{rel780}
\end{eqnarray}
(6) With same notation as in Item (5), we have the following
estimate
\begin{eqnarray}
\|\,T_\varepsilon\otimes {\text{Id}_{\hat{E}}}\,\|_{cb}&\leq &
\sup_\alpha \|\,\tilde{\E_\alpha}(T_\varepsilon\otimes
{\text{Id}_{\hat{E}}})\tilde{\E_\alpha}: L_p(\N_\alpha; \hat{E})
\rightarrow L_p(\N_\alpha; \hat{E})\,\|_{cb}\,.\label{relnice}
\end{eqnarray}
For all $\alpha\geq 1$\,, let
$m_\alpha=\textup{dim}\,(\N_\alpha)$\,. It follows that for $i\in
I$\,,
\begin{equation}\label{rel8798}
L_p(\N_\alpha; E_i)=S_p^{m_\alpha}[\,E_i\,]\,,\quad L_p(\N_\alpha;
\hat{E})=S_p^{m_\alpha}[\,\hat{E}\,]\,.
\end{equation}
Furthermore, by Lemma 5.4 in \cite{P}, there is a complete
isometry
\begin{eqnarray*}
\phi_\alpha: S_p^{m_\alpha}[\,\hat{E}\,]&\rightarrow&
(S_p^{m_\alpha}[E_i])_{\U}\,.
\end{eqnarray*}
By assumption, for all $i\in I$ the operator space $E_i$ is
$OUMD_p$\,. Hence, by (\ref{rel8798}) the map
$\tilde{\E_\alpha}(T_\varepsilon\otimes {\text{Id}_{E_i}})
\tilde{\E_\alpha}: S_p^{m_\alpha}[E_i]\rightarrow
S_p^{m_\alpha}[E_i]$ is completely bounded. Moreover, since
$\tilde{\E_\alpha}$ is a complete contraction, we obtain the
estimate
\begin{eqnarray}
\|\, \tilde{\E_\alpha}(T_\varepsilon\otimes {\text{Id}_{E_i}})
\tilde{\E_\alpha}\,\|_{cb}\,\,\leq &\|\,T_\varepsilon\otimes
{\text{Id}_{E_i}}\,\|_{cb}&\leq\,\,c_p(E_i)\,.\label{rel8900}
\end{eqnarray}
Define a map $\psi_\alpha$ by
\begin{eqnarray*}
\psi_\alpha&=&\left(\,\tilde{\E_\alpha}(T_\varepsilon\otimes
{\text{Id}_{E_i}}) \tilde{\E_\alpha}\,\right)_{\U}\,.
\end{eqnarray*}
By Proposition 10.3.2 in \cite{ER} it follows that
%\begin{eqnarray*}
$\psi_\alpha\in \mathcal{C}\B(\,(S_p^{m_\alpha}[\,E_i\,])_{\U},
(S_p^{m_\alpha}[\,E_i\,])_{\U}\,)$
 and, moreover,
$\|\psi_\alpha\|_{cb}\leq \sup_i \,c_p(E_i)\,.$ Therefore, we
obtain  the following commuting diagram of completely bounded maps
$$
\xymatrix{
%&{\M}&\\
{S_p^{m_\alpha}[\hat{E}]}
\ar@{->}^{\tilde{\E_\alpha}(T_\varepsilon\otimes
{\text{Id}_{\hat{E}}})\tilde{\E_\alpha}}[rr]
\ar@<-1ex>^{\phi_\alpha}[dd]
%\ar@{^{(}->}[ur]
& & S_p^{m_\alpha}[\hat{E}]
% \ar@<-2ex>^{\phi_m^{-1}}
%{{\E}_{{{M_{n_2}}{{\M}_2}}}}
%{\mbox{ \tiny
%cond.}\atop{\mbox{ \tiny expect.}}}}
%[dd]
%\ar@{_{(}->}[ul]
\\
&&\\
(S_p^{m_\alpha}[E_i])_{\U}
% \ar@{^{(}->}[uu] \ar_{*\mbox{\tiny -hom }\Phi_1}[rdd]
\ar^{\psi_\alpha}[rr] & & (S_p^{m_\alpha}[E_i])_{\U}
\ar@<-1ex>^{\phi_\alpha^{-1}}[uu] & }
$$
%\\
%&&\\
%&
%{A_1\subseteq B_1}\atop
%{B_1=A_1\oplus r_1\mathbb{C}}
%\ar@{^{(}->}[uur]_{\mbox{\tiny unital}}
%&
%}
%$$
We deduce that
\begin{equation*}
\|\,\tilde{\E_\alpha}(T_\varepsilon\otimes {\text{Id}_{\hat{E}}})
\tilde{\E_\alpha}\,\|_{cb}\leq  \|\,\phi_\alpha^{-1}\circ
\psi_\alpha\circ \phi_\alpha\,\|_{cb} \leq
\|\,\phi_\alpha^{-1}\,\|_{cb}\|\,\psi_\alpha\,\|_{cb}\|\,\phi_\alpha\,\|_{cb}
\leq  \sup_i c_p(E_i)\,.
\end{equation*}
By (\ref{relnice}), this yields the conclusion. Moreover, the
$OUMD_p$ constant of the ultraproduct $\hat{E}$ satisfies
$c_p(\hat{E})\leq \sup_i c_p(E_i)$\,.
%\begin{eqnarray*}
%c_p(\hat{E})&\leq &\sup_i c_p(E_i)\,.
%\end{eqnarray*}
By general results on ultraproducts, it follows that
\begin{equation*}
c_p(\hat{E})\leq \lim\limits_{\U} c_p(E_i)\,.
\end{equation*}
(7) We first show that $E\oplus_p F$ is $OUMD_p$ with respect to
$\M$ and the filtration $(\M_n)_{n\geq 1}$\,. By Remark
\ref{umdpislocal}, it is enough to prove that there exists $c_p(E,
F)> 0$ so that for all $m\geq 1$\,,
\begin{eqnarray}\label{eq98986}
\|\,T_\varepsilon\otimes {\text{Id}_{E\oplus_p F}}: S_p^m[
E\oplus_p F]\rightarrow S_p^m[E\oplus_p F]\,\|_{cb}&\leq &c_p(E,
F)\,.
\end{eqnarray}
Recall that by (2.9) in \cite{P} we have the complete isometry
\begin{eqnarray*}
S_p^m[E\oplus_p F]=S_p^m[E]\oplus_p {S_p^m[F]}\,.\label{spoplus}
\end{eqnarray*}
It follows that
\begin{eqnarray*}
\|\,T_\varepsilon\otimes {\text{Id}_{E\oplus_p F}}:
S_p^m[E\oplus_p F]\rightarrow S_p^m[E\oplus_p F]\,\|_{cb}&\leq &
{[(c_p(E))^p+(c_p(F))^p]}^{\frac1p}\,,
\end{eqnarray*}
where $c_p(E)$ and $c_p(F)$ are the $OUMD_p$ constants of $E$\,,
respectively $F$\,. Hence (\ref{eq98986}) is proved. Furthermore,
note that if $1< q< \infty$\,, then
\begin{eqnarray*}
d_{cb}(E\oplus_p F\,, E\oplus_q F)&\leq &2\,.
%\label{dcboplus}
\end{eqnarray*}
The assertion follows now from Item (1). Moreover, from the proof
we obtain the estimate \[ c_p(E\oplus_q F)\leq
{[(c_p(E))^p+(c_p(F))^p]}^{\frac1p}\,. \] (8) The statement
follows immediately by interpolation from Remark
\ref{charactumdp}, using the following completely isometric
identity (see (3.5) in \cite{P})
\begin{eqnarray*}
{[L_q(\M; E_0), L_s(\M; E_1)]}_\theta&=&L_p(\M; E_\theta)\,.
\end{eqnarray*}
Moreover, the corresponding constants satisfy the estimate
\begin{eqnarray}
c_p(E_\theta)&\leq
&{[c_q(E_0)]^{1-\theta}}{[c_s(E_1)]^\theta}\,.\label{rel779}
\end{eqnarray}
This completes the proof.\qed
%\\[-1.4cm]\qed

\begin{exam}\label{mnumdp}\rm
If $1< p< \infty\,,$ then the matrix algebras $M_m\,$ are
$OUMD_p$\,, for all positive integers $m\,.$ This follows
immediately from Remark \ref{Complumdp} and (4) in Proposition
\ref{propertiesumdp}.
\end{exam}

\vspace*{0.3cm} The following lemma provides a necessary condition
for an operator space $E$ to be $OUMD_p\,.$
\begin{lemma}\label{connectumd}
Let $1< p< \infty\,,$ and $E$ an operator space. If $E$ is
$OUMD_p\,,$ then $S_p[E]$ is $UMD$ (as a Banach space).
\end{lemma}
 \vspace*{0.3cm}\noindent {\em Proof.} Since $UMD$ is a local
property, it suffices to show that for all positive integers
$m\,,$ the Banach space $S_p^m[E]$ is $UMD$. We will prove that
there exists a constant $k_p> 0\,,$ independent of $m$\,, such
that
\begin{eqnarray}
\|\,T_{\varepsilon}\otimes {\text{Id}_{S_p^m[E]}}: L_p(([0, 1],
{\F}, \mu); S_p^m[E])\rightarrow L_p(([0, 1], {\F}, \mu);
S_p^m[E])\,\|&\leq & k_p\,.\label{relat101}
\end{eqnarray}
Here $\mu$ is the Lebesgue measure on $[0, 1]\,, {\F}$ a
$\sigma$-algebra of subsets of $[0, 1]\,,$ and $({\F}_n)_{n\geq
1}$ a filtration of $\F\,.$ Martingale differences are considered
with respect to this filtration. By Fubini's theorem, we have the
complete isometry
\begin{eqnarray}
L_p([0, 1]; S_p^m[E])&=&L_p([0, 1]; L_p(M_m; E))\\
\nonumber
&=&L_p(L_\infty([0, 1]); L_p(M_m; E))\\
\nonumber &=&L_p(L_\infty([0, 1])\otimes {M_m}; E)\label{fubistep}
\end{eqnarray}
Note that $\,(L_\infty([0, 1], {\F}_n, \mu)\otimes {M_m})_{n\geq
1}$ is a filtration of the hyperfinite von Neumann algebra
$L_\infty([0, 1])\otimes M_m$\,. Moreover, for all positive
integers $n$\,, the unique trace-preserving conditional
expectation onto the subalgebra $L_\infty([0, 1], {\F}_n)\otimes
M_m$ is $\,{\E}_n=\mathbb{E}_n\otimes {\text{Id}_{M_m}}\,,$ where
$\mathbb{E}_n=\mathbb{E}(\,\cdot\, |{\F}_n)\,.$   Let
$k_p=c_p(E)$\,, where $c_p(E)$ is the $OUMD_p$ constant of $E$\,.
The fact that $E$ is $OUMD_p$ with respect to the filtration
$(L_\infty([0, 1], {\F}_n, \mu)\otimes {M_m})_{n\geq 1}$ of
$L_\infty([0, 1])\otimes M_m$\,, together with (2.19) yields
(\ref{relat101}), and the proof is complete. \qed

 \noindent
\begin{exam}\rm
There exists a Hilbert space, subspace of some commutative
$C^*$-algebra which is $UMD$ (as a Banach space), but not
$OUMD_p$\,, for any $1< p< \infty\,.$
\end{exam}
\vspace*{0.3cm} Let $E=\min(l_2)\,.$ Then $E$ is a Hilbert space,
hence it is $UMD$ (as a Banach space).
%We will show that $E$ is not $OUMD_p$, for any $\,1< p< \infty\,.$
Assume that $E$ is $OUMD_p\,,$ for some $1< p< \infty\,.$ We will
show that this leads to a contradiction. Following Pisier's
argument in \cite{P}, let $(\M,\tau)$ be the hyperfinite $II_1$
factor equipped with the canonical "dyadic filtration"
$({\M}_n)_{n\geq 1}$, where ${\M}_n={M_2\otimes\ldots\otimes M_2}$
($n$ times). By (3) in Proposition \ref{propertiesumdp}, it
follows that $L_p(\M; E)$ is $OUMD_p$\,; hence it is $UMD$ (as a
Banach space). However, by Proposition 4.3 in \cite{P}, $L_p(\M;
E)$ contains a subspace isomorphic to $c_0\,.$ This contradicts
the fact that $L_p(\M; E)$ is $UMD$ (as a Banach space), since
$c_0$ is not $UMD$ (being non-reflexive). Therefore $E$ is not
$OUMD_p$\,, for any $1< p< \infty\,.$

\section{Main results}

As a consequence of Proposition \ref{tepscb}, together with
Fubini's theorem, we immediately obtain the following

\begin{prop}\label{lpumdp}
Let $1< p< \infty\,.$ Then $L_p(\N, \phi)$ is $OUMD_p$\,, for
every von Neumann algebra $\N$ equipped with an n.s.f. tracial
state $\phi$\,.
\end{prop}

\re \label{spumdp} Let $1< p< \infty$ and $m$ be a positive
integer. It follows that $S_p=L_p(\B(l_2), \text{tr})$ and,
respectively, $S_p^m=L_p(M_m, \text{Tr}_m)$ are $OUMD_p$\,. Here
$\text{tr}$ denotes the usual trace on $\B(l_2)$\,, while
$\text{Tr}_m$ denotes the non-normalized trace on $M_m\,.$
Moreover, the corresponding $OUMD_p$ constants satisfy $\sup_m
c_p(S_p^m)\,\leq\, c_p(S_p)\,< \,\infty\,.$\mar

\begin{prop}\label{ohumdp}
The operator Hilbert space $OH$ is $OUMD_p$\,, for $1< p<
\infty\,.$
\end{prop}
\vspace*{0.3cm}
 \noindent {\em Proof.} Recall that by Corollary
2.6 in \cite{Pi2}\,, $OH$ can be obtained by complex interpolation
between the column Hilbert space $C$ and the row Hilbert space
$R$, namely,
\begin{eqnarray}
OH&=&[C, R]_{\frac12} \,.
\end{eqnarray}
By the reiteration theorem for the complex method (\cite{BL},
Theorem 4.6.1), it follows from (\ref{cpandcpm}) and
(\ref{rpandrpm}) that for $1< p< \infty$ we have
\begin{eqnarray}
OH&=&[C_p, R_p]_{\frac12}\,.\label{ohhalf}
\end{eqnarray}
Note that we can view $C_p$ as the column space of $S_p$ and,
respectively, we can view $R_p$ as the row space of $S_p\,.$
Hence, by (2) in Proposition \ref{propertiesumdp}, and Remark
\ref{spumdp} it follows that both $C_p$ and $R_p$ are $OUMD_p$\,.
Therefore, by (\ref{ohhalf}), a further application of (8) in
Proposition \ref{propertiesumdp} yields the conclusion.\qed

\begin{cor}\label{findimumdp}
Every finite dimensional operator space $E$ is $OUMD_p$\,, for $1<
p< \infty\,.$
\end{cor}
\vspace*{0.3cm} \noindent {\em Proof.} Let $E$ be an
$n$-dimensional operator space. By Corollary 9.3 of \cite{Pi2},
there is an isomorphism $\phi: {OH}_n\rightarrow E\,,$ such that
$\|\phi\|_{cb}\|\phi^{-1}\|_{cb}\leq \sqrt{n}\,.$ Thus the
c.b.Banach-Mazur distance between $E$ and ${OH}_n$ satisfies
\[ d_{cb}( E, {OH}_n)\,\leq \,\sqrt{n}\,. \]
Hence $E$ is completely isomorphic to ${OH}_n\,.$ The conclusion
follows by Proposition \ref{ohumdp} (applied in the finite
dimensional case) and (1) in Proposition \ref{propertiesumdp}.\qed

\begin{theorem}\label{cuumdp}
If $1< u, p< \infty\,,$ then $C_u$ is $OUMD_p\,.$
\end{theorem}
\vspace*{0.3cm} \noindent {\em Proof.} Since $OUMD_p$ is a local
property (see Remark \ref{umdpislocal}),
%by the density of  $(C_u^m)_{m\geq 1}$ in $C_u$\,,
it suffices to show that $C_u^m$ is $OUMD_p$\,, for all positive
integers $m$, and
\[ \sup_m c_p(C_u^m)< \infty\,. \] Let $\,(\M,
\tau)\,$ be a hyperfinite von Neumann algebra with an n.f. tracial
state $\tau$ and $(\M_n)_{n\geq 1}$ a filtration of $\M$\,.
%We will show that $C_u^m$ is $OUMD_p$ with respect to $(\M, \tau)$ and the filtration $(\M_n)_{n\geq 1}\,.$\\
By Remark \ref{spumdp}, $S_p^m$ is $OUMD_p$ with respect to the
filtration $(\M_n)_{n\geq 1}$ of $\M$. By (2) in Proposition
\ref{propertiesumdp}, the same holds for $C_p^m$ and $R_p^m$\,, as
subspaces of $S_p^m\,.$ We claim that there exist $0< \theta< 1$
and $1< q, s< \infty$ such that the following complete isometry
holds
\begin{eqnarray}
{[L_q(\M; C_q^m), L_s(\M; R_s^m)]}_\theta&=&L_p(\M;
C_u^m)\,.\label{lqspinterp}
\end{eqnarray}
By applying (3.5) in \cite{P}, this reduces to showing the
existence of $0< \theta< 1$ and $1< q, s< \infty$ such that the
following two relations hold
\begin{eqnarray}
\frac1p&=&{\frac{1-\theta}{q}}+{\frac{\theta}{s}}\label{pqr}\,,\\
C_u^m&=&{[C_q^m, R_s^m]}_\theta\,.\label{cqrscu}
\end{eqnarray}
Let $s'$ denote the conjugate exponent of $s$\,. Then we have
\[ R_s^m\,\,=\,\,[R^m, C^m]_{\frac1s}\,\,=\,\,[C^m, R^m]_{1-{\frac1s}}\,\,=\,\,C_{s'}^m\,, \]
Therefore, relations (\ref{pqr}) and (\ref{cqrscu}) are equivalent
to
\[ \left\{\begin{array}{lcl}
                                 \frac1p&=&{\frac{1-\theta}{q}}+{\frac{\theta}{s}}\,,\vspace{0.3cm}\\
                                 \frac1u&=&{\frac{1-\theta}{q}}+{\frac{\theta}{s'}}\,.\\
                                \end{array}
                        \right. \]
Equivalently,
\begin{eqnarray*}
{\frac1p}+{\frac1u}-\theta=\frac{2(1-\theta)}{q}\,\,&\mbox{and}&
\,\,{\frac1p}-{\frac1u}+\theta=\frac{2{\theta}}{s}\,.
\end{eqnarray*}
Thus, we have to show that there exists $\,0< \theta< 1$ such that
\[ {\frac1p}+{\frac1u}-\theta\,>0\,,\quad
2(1-\theta)\,>{\frac1p}+{\frac1u}-\theta\,,\quad
{\frac1p}-{\frac1u}+\theta\,> 0\,,\quad
2\,{\theta}\,>{\frac1p}-{\frac1u}+\theta\,.\]
%Equivalently, we
%have to show that there exists $\,0< \theta< 1$ such that all four
%conditions
%\[ \theta\,< {\frac1p}+{\frac1u}\,, \quad \theta\,< 2-({\frac1p}+{\frac1u})\,, \quad \theta>{\frac1u}-{\frac1p}\,, \quad \theta\,>{\frac1p}-{\frac1u}\,.\]
%are satisfied.
These conditions are equivalent to
\begin{eqnarray}
\left|\,{\frac1p}-{\frac1u}\,\right|\,< &\theta&<
\,\min\left\{{\frac1p}+{\frac1u}\,,
2-\left({\frac1p}+{\frac1u}\right)\right\}\,.\label{existteta}
\end{eqnarray}
Since $\,1< p, u< \infty$\,, the following relations hold
\begin{eqnarray}
\left|\,{\frac1p}-{\frac1u}\,\right|\,\,<\,\, \max\left\{\frac1p\,, \frac1u\right\}&<& 1\,,\\
\min\left\{\,{\frac1p}+{\frac1u}\,,
2-\left({\frac1p}+{\frac1u}\right)\,\right\}&>& 0\,.
\end{eqnarray}
Since $\,1< p, u$ we also have
\begin{eqnarray*}
{\frac1p}-{\frac1u}<
2-\left(\,{\frac1p}+{\frac1u}\,\right)\,\,&\mbox{and}&\,\,{\frac1u}-{\frac1p}<
2-\left(\,{\frac1p}+{\frac1u}\,\right)\,.
\end{eqnarray*}
Therefore, we see that
\begin{eqnarray*}
0\,\,\,<\,\,\, \left|\,
{\frac1p}-{\frac1u}\,\right|&<&\min\left\{{\frac1p}+{\frac1u}\,,
2-\left({\frac1p}+{\frac1u}\right)\,\right\}\,,
\end{eqnarray*}
which implies the existence of some $0< \theta< 1$ satisfying
(\ref{existteta})\,. Further, set
\begin{eqnarray*}
q\,=\,\frac{2(1-\theta)}{{\frac1p}+{\frac1u}-\theta}\,\,&\mbox{and}&\,\,s\,=\,\frac{2\,\theta}{{\frac1p}-{\frac1u}+\theta}\,.
\end{eqnarray*}
The above relations ensure that $1< q, s< \infty\,$ and therefore
the claim is proved. Since $C_q^m$ is $OUMD_q$ and $R_s^m$ is
$OUMD_s$\,, it follows by interpolation from (\ref{lqspinterp})
that $C_u^m$ is $OUMD_p$ with respect to the filtration
$(\M_n)_{n\geq 1}$ of $(\M, \tau)$\,. This argument implies that
$C_u^m$ is $OUMD_p$\,. Moreover, using (\ref{rel779}) and
(\ref{rel778}), we obtain from the proof the following estimates
for the corresponding $OUMD_p$ constants
\begin{equation*}
c_p(C_u^m)\leq [c_q(C_q^m)]^{1-\theta} [c_s(R_s^m)]^\theta \leq
[c_q(S_q)]^{1-\theta} [c_s(S_s)]^\theta\,. \end{equation*} This
implies that $\,\sup_m c_p(C_u^m)< \infty\,,$ and the conclusion
follows.\qed

\re\label{umdru} Let $1< u, p< \infty\,.$ By (\ref{cpandcpm}) it
follows from the equivalence theorem for the complex method
(Theorem 4.3.1 in \cite{BL}) that
\begin{eqnarray}
(C_u)^*\,\,=\,\,R_{u'}&\mbox{and}&(R_u)^*\,\,=\,\,C_{u'}\,,\label{dualcuru}
\end{eqnarray}
where ${\frac1u}+{\frac1{u'}}=1\,.$ Therefore, by Theorem
\ref{cuumdp} and (5) in Proposition \ref{propertiesumdp} it
follows that $R_u$ is $OUMD_p\,.$ \mar

\begin{prop}\label{colrowhaag1}
Let $1< u, p< \infty\,.$ Then the spaces $C_p\,{\otimes^{h}}
C_u$\,, $C_u\,{\otimes^{h}}R_p$\,, $C_p\,{\otimes^{h}}R_u$\,,
$R_u\,{\otimes^{h}}R_p$\, are all $OUMD_p$\,.
\end{prop}
 \vspace*{0.3cm}\noindent {\em Proof.} Let $m$ be a positive
integer, $(\M, \tau)$ a hyperfinite von Neumann algebra equipped
with an n.s.f. tracial state $\tau$ and $(\M_n)_{n\geq 1}$ a
filtration of $\M$\,. By Theorem \ref{cuumdp} it follows that, in
particular, $C_u^m$ is $OUMD_p$ with respect to the filtration
$(\M_n\otimes M_m)_{n\geq 1}$ of $(\M\otimes M_m, \tau\otimes
{\text{tr}_m})\,.$ Therefore,
\begin{equation}\label{eq987678}
\!\!\!\!\|\,T_\varepsilon\otimes {\text{Id}_{C_u^m}}:
L_p(\M\otimes M_m; C_u^m)\rightarrow L_p(\M\otimes M_m;
C_u^m)\,\|_{cb}\leq c_p(C_u^m)\,\,\leq \,\,c_p(C_u)\,\,<\,\,
\infty\,.
\end{equation}
Fubini's theorem and (\ref{spehaage}) yield the complete
isometries
\begin{eqnarray*}
L_p(\M\otimes M_m; C_u^m)&=&L_p(\M; L_p(M_m; C_u^m))\\
&=&L_p(\M; S_p^m[C_u^m])\\
&=&L_p(\M; C_p^m\,{\otimes^{h}}C_u^m\,{\otimes^{h}}R_p^m)\,.
\end{eqnarray*}
By the injectivity of the Haagerup tensor product it follows that
\begin{eqnarray*}
L_p(\M; C_p^m\,{\otimes^{h}}C_u^m)&\subseteq&L_p(\M; C_p^m\,{\otimes^{h}}C_u^m\,{\otimes^{h}}R_p^m)\,;\\
L_p(\M; C_u^m\,{\otimes^{h}}R_p^m)&\subseteq&L_p(\M;
C_p^m\,{\otimes^{h}}C_u^m\,{\otimes^{h}}R_p^m)\,.
\end{eqnarray*}
By (\ref{eq987678}), it follows that $C_p^m\,{\otimes^{h}}C_u^m$
and $C_u^m\,{\otimes^{h}}R_p^m$ are both $OUMD_p$ with respect to
the filtration $(\M_n)_{n\geq 1}$ of $(\M, \tau)\,.$ Since this
filtration was arbitrarily chosen, this argument shows that
$C_p^m\,{\otimes^{h}}C_u^m$ and $C_u^m\,{\otimes^{h}}R_p^m$ are
both $OUMD_p$\,. Moreover, from the proof we obtain the following
estimates for the corresponding constants
\begin{eqnarray}
\!\!\!\!\!\!c_p(C_p^m\,{\otimes^{h}}C_u^m)\,,
\,\,c_p(C_u^m\,{\otimes^{h}}R_p^m)&\leq
&c_p(C_u^m)\,\,\leq\,\,c_p(C_u)\,\,<\,\,\infty\,.\label{rel210}
\end{eqnarray}
%By density, this implies that the spaces $\,C_p\,{\otimes^{h}}C_u\,$ and $\,C_u\,{\otimes^{h}}R_p\,$ are $OUMD_p$ with respect %to $(\M, \tau)$ and the filtration $(\M_n)_{n\geq 1}\,.$
A similar argument applied to the space $R_u^m$, which is $OUMD_p$
as a subspace of $R_u$ (see Remark \ref{umdru}), shows that
$C_p^m\,{\otimes^{h}}R_u^m$ and $R_u^m\,{\otimes^{h}}R_p^m$ are
both $OUMD_p\,.$ Moreover,
\begin{eqnarray}
\!\!\!\!\!\!c_p(C_p^m\,{\otimes^{h}}R_u^m)\,,
\,\,c_p(R_u^m\,{\otimes^{h}}R_p^m)&\leq
&c_p(R_u^m)\,\,\leq\,\,c_p(R_u)\,\,<\,\,\infty\,.\label{rel222}
\end{eqnarray}
Since $m$ is arbitrarily chosen, the conclusion follows by density
and the fact that $OUMD_p$ is a local property.\qed

\begin{prop}\label{colrowhaag2}
Let $1< u, p< \infty\,.$ Then the spaces
$R_p\,{\otimes^{h}}C_u$\,, $C_u\,{\otimes^{h}}C_p$\,,
$R_p\,{\otimes^{h}}R_u$\, and $R_u\,{\otimes^{h}}C_p\,$ are all
$OUMD_p$\,.
\end{prop}
 \vspace*{0.3cm}\noindent {\em Proof.} Let $m$ be a positive
integer. By Proposition \ref{colrowhaag1}, the operator spaces
$C_{p'}^m\otimes^{{h}} C_{u'}^m$\,, $C_{u'}^m\otimes^{{h}}
R_{p'}^m$\,, $C_{p'}^m\otimes^{{h}} R_{u'}^m$ and
$R_{u'}^m\otimes^{{h}} R_{p'}^m$ are all $OUMD_{p'}$\,, where
${\frac1p}+{\frac1{p'}}=1$ and ${\frac1u}+{\frac1{u'}}=1$\,.
Hence, by (5) in Proposition \ref{propertiesumdp}, the dual spaces
are $OUMD_p$\,. By (\ref{dualcuru}) and the self-duality of the
Haagerup tensor product in the finite dimensional case (see
\cite{ER1}), we obtain the complete isometries
\begin{eqnarray*}
(C_{p'}^m\otimes^{{h}} C_{u'}^m)^*\,\,=\,\,R_{p}^m\otimes^{{h}}
R_u^m\,, && (C_{u'}^m\otimes^{{h}}
R_{p'}^m)^*\,\,=\,\,R_u^m\otimes^{{h}} C_p^m\,,
\end{eqnarray*}
respectively,
\begin{eqnarray*}
(C_{p'}^m\otimes^{{h}} R_{u'}^m)^*\,\,=\,\,R_p^m\otimes^{{h}}
R_u^m\,, && (R_{u'}^m\otimes^{{h}}
R_{p'}^m)^*\,\,=\,\,C_u^m\otimes^{{h}} C_p^m\,.
\end{eqnarray*}
Moreover, from (\ref{rel210}), (\ref{rel222}) and (\ref{rel780})
we obtain the following estimates for the corresponding $OUMD_p$
constants
\begin{eqnarray*}
\!\!\!\!\!\!c_p(R_p^m\,{\otimes^{h}}R_u^m)\,, \,\,c_p(R_u^m\,{\otimes^{h}}C_p^m)&\leq &c_{p'}(C_u^m)\,\,\leq\,\,c_{p'}(C_u)\,\,<\,\,\infty\,,\\
\!\!\!\!\!\!c_p(R_p^m\,{\otimes^{h}}R_u^m)\,,
\,\,c_p(C_u^m\,{\otimes^{h}}C_p^m)&\leq
&c_{p'}(R_u^m)\,\,\leq\,\,c_{p'}(R_u)\,\,<\,\,\infty\,.
\end{eqnarray*}
The conclusion follows from the fact that $OUMD_p$ is a local
property.\qed

\begin{prop}\label{spumdintermed}
Let $1< p< \infty\,.$ If $1< q< \infty$ satisfies $\frac{2p}{p+1}<
q< 2p\,,$ then $S_q$ is $OUMD_p$\,.
\end{prop}
\vspace*{0.3cm} \noindent {\em Proof.} We first show that there
exist $0< \theta< 1$ and $1< u, v< \infty$ such that
\begin{eqnarray}
{[C_u\,, C_p]}_\theta=C_q&\mbox{and}&{[R_p\,,
R_v]}_\theta=R_q\,.\label{findteta}
\end{eqnarray}
This is equivalent to showing that there exist $0< \theta< 1$ and
$1< u, v< \infty$ such that the following relations hold
%\[ \left\{\begin{array}{lcl}
                                % \frac1q&=&{\frac{1-\theta}{u}}+{\frac{\theta}{p}}\vspace{0.3cm}\\
                                 %\frac1q&=&{\frac{1-\theta}{p}}+{\frac{\theta}{v}}\\
                                %\end{array}
                        %\right. \]
\begin{eqnarray}
\frac1q&=&{\frac{1-\theta}{u}}+{\frac{\theta}{p}}\,,\label{rel11}\\
 \frac1q&=&{\frac{1-\theta}{p}}+{\frac{\theta}{v}}\,.\label{rel12}
\end{eqnarray}
By the assumption on $q$ we have $\frac{2}{q}< 1+{\frac1p}$\,.  An
easy computation shows that
\begin{eqnarray*}
{p'}\left({\frac1q}-{\frac1p}\right)&<& \frac{p'}{q'}\,,
\end{eqnarray*}
where ${\frac1p}+{\frac1{p'}}=1\,$ and
$\,{\frac1q}+{\frac1{q'}}=1\,.$ This ensures the existence of
some $0< \theta< 1$ such that
\begin{eqnarray}
\Big\vert{p'}\left({\frac1q}-{\frac1p}\right)\Big\vert<&\theta&<
\min\left\{1, \frac{p'}{q'}\right\}\,.\label{condteta}
\end{eqnarray}
By (\ref{condteta}) it follows that ${\frac1q}-{\frac{\theta}{p}}<
1-\theta\,.$ This implies the existence of some $1< u< \infty$
such that (\ref{rel11}) holds. Furthermore, (\ref{condteta}) also
shows that $\frac1q< {\frac{1-\theta}{p}}+\theta\,.$ This implies
the existence of some $1< v< \infty$ such that (\ref{rel12})
holds. Therefore the claim is proved.\\ An application of Kouba's
interpolation result (Theorem \ref{koubast}) yields the complete
isometry
\begin{eqnarray}\label{eq345656}
{[C_u\,{\otimes^{h}}R_p\,,
C_p\,{\otimes^{h}}R_v]}_\theta&=&[C_u\,,
C_p]_\theta\,{\otimes^{h}}[R_p\,, R_v]_\theta\\\nonumber
&=&C_q\,{\otimes^{h}}R_q\,\,=\,\,S_q\,.
\end{eqnarray}
By Proposition \ref{colrowhaag1}, both spaces
$C_u\,{\otimes^{h}}R_p$ and  $C_p\,{\otimes^{h}}R_v$ are
$OUMD_p$\,. Therefore, by (8) in Proposition \ref{propertiesumdp},
it follows by interpolation from (\ref{eq345656}) that $S_q$ is
$OUMD_p$\,. \qed

\begin{prop}\label{spumdpprime}
If $1< p< \infty$ and ${\frac1p}+{\frac1{p'}}=1$\,, then $S_{p'}$
is $OUMD_p$\,.
\end{prop}
\vspace*{0.3cm} \noindent {\em Proof.}
%By Kouba's interpolation
%result (Theorem \ref{koubast}), if $\,0< \theta< 1$ and $\,1< u,
%v< \infty$, we have the complete isometry
%\begin{eqnarray}
%{[C_u\,{\otimes^{h}}R_p\,,
%C_p\,{\otimes^{h}}R_v]}_\theta&=&[C_u\,,
%C_p]_\theta\,{\otimes^{h}}[R_p\,, R_v]_\theta\,.
%\end{eqnarray}
{\em Case 1}: $2\leq p$\,. We claim that there exist $0< \theta<
1$ and $1< u, v< \infty$ such that
\begin{eqnarray}
 {[C_u\,, C_p]}_\theta=C_{p'}&\mbox{and}&{[R_p\,, R_v]}_\theta=R_{p'}\,.\label{findtetaprime}
\end{eqnarray}
Equivalently,
 %This is equivalent to showing that there exist $\,0< \theta< 1$ and $\,1< u, v< \infty$ such that
\begin{eqnarray}
\frac1{p'}&=&{\frac{1-\theta}{u}}+{\frac{\theta}{p}}\,,\label{rel3}\\
\frac1{p'}&=&{\frac{1-\theta}{p}}+{\frac{\theta}{v}}\label{rel4}\,.
\end{eqnarray}
Note that (\ref{rel3}) is equivalent to
\begin{eqnarray}\label{rel5}
\frac1p={\frac{1-\theta}{u'}}+{\frac{\theta}{p'}}
&\Leftrightarrow&\frac{1+\theta}{p}={\frac{1-\theta}{u'}}+\theta\,,
\end{eqnarray}
where ${\frac1u}+{\frac1{u'}}=1$\,. This yields a restriction upon
$\theta$ as follows
\begin{eqnarray}
\theta< \frac{1+\theta}{p}< 1&\Leftrightarrow& \frac1{p-1}<
\theta< p-1\,.
\end{eqnarray}
Note that since $2< p< \infty\,,$ we have $0< \frac1{p-1}< 1<
p-1\,.$ Therefore, the above relations yield the restriction
\begin{eqnarray}
\frac1{p-1}<&\theta&< 1\,.\label{restriction1}
\end{eqnarray}
Choose $0< \theta< 1$ such that (\ref{restriction1}) holds. Then
we find $1< u'< \infty$ by solving (\ref{rel5}), and let $u$ be
the conjugate exponent of $u'$\,. Furthermore, note that
(\ref{rel4}) is equivalent to
\begin{eqnarray}
1-{\frac1p}={\frac{1-\theta}{p}}+{\frac{\theta}{v}}&\Leftrightarrow&1-{\frac{2-\theta}{p}}=\frac{\theta}{v}\,.\label{rel6}
\end{eqnarray}
This implies that
\begin{eqnarray}
0< \frac{2-\theta}{p}< 1&\Leftrightarrow&p>2-\theta\,,
\end{eqnarray}
which is obviously true, since by assumption $p> 2> 2-\theta\,.$
Therefore, we can solve (\ref{rel6}) to find $v$\,. Our claim is
completely proved. By Kouba's interpolation result and
(\ref{findtetaprime}) we get
\begin{eqnarray}
{[C_u\,{\otimes^{h}}R_p\,,
C_p\,{\otimes^{h}}R_v]}&=&C_{p'}\,{\otimes^{h}}R_{p'}\,\,=\,\,S_{p'}\,.\label{haagspprime}
\end{eqnarray}
By Proposition \ref{colrowhaag1}, both spaces
$C_u\,{\otimes^{h}}R_p$ and  $C_p\,{\otimes^{h}}R_v$ are
$OUMD_p$\,.
Hence, by (8) in Proposition \ref{propertiesumdp}, it follows by interpolation that $S_{p'}$ is $OUMD_p$\,. \\
{\em Case 2}: $1< p< 2\,.$ The assertion follows by duality.
Indeed, since $2< p'< \infty\,,$ it follows by Case 1 that
$S_{(p')'}=S_p$ is $OUMD_{p'}\,.$  By (5) in Proposition
\ref{charactumdp}, this implies that $S_{p'}$ is $OUMD_p\,.$\qed

\begin{lemma}\label{squmd2}
If $1< q< \infty\,,$ then $S_q$ is $OUMD_2$\,.
\end{lemma}
\vspace*{0.3cm} \noindent {\em Proof.} By Remark \ref{spumdp},
$S_q$ is $OUMD_q$\,. Also, by Proposition \ref{spumdpprime}, $S_q$
is $OUMD_{q'}$\,, where ${\frac1q}+{\frac1{q'}}=1$. By (8) in
Proposition \ref{propertiesumdp}, interpolation with exponent
$\theta=\frac12$ yields the conclusion.\qed

\begin{theorem}\label{finalsqumdp}
If $1< p, q< \infty\,,$ then $S_q$ is $OUMD_p$\,.
\end{theorem}
\vspace*{0.3cm} \noindent {\em Proof.} {\em Case 1}: $2\leq q<
\infty\,.$ We will first show that if $p\geq 2\,,$ then $S_q$ is
$OUMD_p$\,. Indeed, if $p\geq q$\,, then we have $\frac{2p}{p+1}<
2\leq q< 2p\,.$ Thus, by Proposition \ref{spumdintermed}, it
follows that $S_q$ is $OUMD_p\,.$ The case $p=2$ follows from
Lemma \ref{squmd2}. On the other hand, if $2< p< q$\,, then there
exists $0< \theta< 1$ such that
\begin{eqnarray}
\frac1p&=&{\frac{1-\theta}{2}}+{\frac{\theta}{q}}\,.\label{relp2q}
\end{eqnarray}
By  Lemma \ref{squmd2}, $S_q$ is $OUMD_2\,.$ Also, by Remark
\ref{spumdp}, $S_q$ is $OUMD_q\,.$ Interpolating with exponent
$\theta$ given by (\ref{relp2q})\,, it follows from (8) in
Proposition \ref{propertiesumdp} that $S_q$ is $OUMD_p\,,$ which
proves the claim. In particular, as a consequence it follows by
duality that $S_2$ is $OUMD_u$\,, for all $1< u< \infty$\,. Hence
the case $q=2$ is completely proved.\\
Furthermore, note that for $2< q< \infty$ we have $1< {q'}< 2\,,$
where $\frac1q+{\frac1{q'}}=1\,.$ Assume that $q'< p< 2\,$. There
exists $0< \eta< 1$ such that
\begin{eqnarray}
\frac1p&=&{\frac{1-\eta}{q'}}+{\frac{\eta}{2}}\,.\label{rel341}
\end{eqnarray}
By Proposition \ref{spumdpprime}, $S_q$ is $OUMD_{q'}\,.$ Also, by
Lemma \ref{squmd2}, $S_q$ is $OUMD_2$\,. Thus, by (8) in
Proposition \ref{propertiesumdp}, interpolation with exponent
$\theta$ given by (\ref{rel341}) shows that $S_q$ is $OUMD_p\,.$
To summarize, we have proved so far that if $p\geq {q'}\,,$ then
$S_q$ is $OUMD_p\,.$ It remains to analyze the case when $1< p<
{q'}\,.$ In this case we have ${q'}< 2< {p'}\,,$ where
$\frac1p+{\frac1{p'}}=1\,.$ By what we proved above, it follows
that $S_q$ is $OUMD_{p'}\,.$ By Proposition \ref{spumdpprime},
$S_p$ is $OUMD_{p'}\,.$ Note that $p< {q'}< q\,,$ hence there
exists $0< \xi< 1$ such that
\begin{eqnarray}
\frac1{q'}&=&{\frac{1-\xi}{p}}+{\frac{\xi}{q}}\,.\label{rel342}
\end{eqnarray}
By (8) in Proposition \ref{propertiesumdp}, interpolation with exponent $\theta$ given by (\ref{rel342}) implies that $S_{q'}$ is $OUMD_{p'}\,.$ By duality, an application of (5) in Proposition \ref{propertiesumdp} shows that $S_q$ is $OUMD_p\,.$ \\
{\em Case 2}: $1< q< 2\,.$ Note that $2\leq q'< \infty\,.$ Thus,
for $1< p< \infty\,,$ it follows by Case 1 that $S_{q'}$ is
$OUMD_{p'}\,,$ where $p'$ is the conjugate exponent of $p$\,. By
duality, (5) in Proposition \ref{propertiesumdp} implies that
$S_q$ is $OUMD_p\,.$ This completes the proof.\qed

\begin{prop}\label{sqsuumdp}
If $1< p, q, u< \infty\,,$ then $S_q[S_u]$ is $OUMD_p$\,.
\end{prop}
\vspace*{0.3cm} \noindent {\em Proof.} We first show that
$S_q[S_u]$ is $OUMD_{u'}$\,, where ${\frac1u}+ {\frac1{u'}}=1\,.$
As showed in the proof of Theorem \ref{cuumdp}, there exist $0<
\theta< 1$ and $1< v, s< \infty$\,, such that the following
relations hold
\begin{eqnarray}
\frac1q&=&{\frac{1-\theta}{v}}+{\frac{\theta}{s}}\,,\label{rel6767}\\
\frac1u&=&{\frac{1-\theta}{v}}+{\frac{\theta}{s'}}\,.\label{rel6768}
\end{eqnarray}
Here $s'$ denotes the conjugate exponent of $s$\,. It follows that
$S_q={[S_v, S_s]_\theta}$ and, respectively $S_u={[S_v,
S_{s'}]_\theta}$\,. An application of Corollary 1.4 in \cite {P}
shows that
\begin{eqnarray}
S_q[S_u]&=&{[\,S_v[S_v]\,, S_s[S_{s'}]\,]}_\theta\,.\label{86}
\end{eqnarray}
Let $v'$ be the conjugate exponent of $v\,.$ By Proposition
\ref{finalsqumdp} it follows that $S_v[S_v]\simeq
S_v(\mathbb{N}\times \mathbb{N})$ is $OUMD_{v'}$\,. Also, since
$S_{s'}$ is $OUMD_s$\,, it follows by (3) in Proposition
\ref{propertiesumdp} that for all $m\geq 1$\,, $S_s^m[S_{s'}]$ is
$OUMD_s\,,$ and moreover, the corresponding constants satisfy
$c_s(S_s^m[S_{s'}])\leq c_s(S_{s'})$\,. Since $OUMD_{s'}$ is a
local property, this implies that $S_s[S_{s'}]$ is $OUMD_s\,.$
Furthermore, from (\ref{rel6768}) we easily deduce that
\begin{eqnarray}\label{89}
\frac1{u'}&=&{\frac{1-\theta}{v'}}+{\frac{\theta}{s}}\,.
\end{eqnarray}
Therefore, by (\ref{86}) and (\ref{89}) it follows by
interpolation that $S_q[S_u]$ is $OUMD_{u'}$\,. Using the fact
that $S_u$ is $OUMD_q$\,, a similar argument as above shows that
$S_q[S_u]$ is $OUMD_q\,.$ By interpolation it follows that $S_q[S_u]$ is $OUMD_p$\,, for $\min\{u', q\}\leq p\leq \max\{u',q\}\,.$ \\
Next, assume that $\max\{u', q, 2\}< p< \infty\,.$ This implies
that $p> \min\{2, u, q\}\,.$ Let $1< t_0< \min\{2, u, q\}< p<
{t_1}< \infty$\,. There exists $0< \eta< 1$ such that
\begin{eqnarray}
\frac1q&=&{\frac{1-\eta}{{t_1}}}+{\frac{\eta}{t_0}}\,.
\end{eqnarray}
Let $\,2< w< \infty$ be such that
\begin{eqnarray}
\frac1u&=&{\frac{1-\eta}{w}}+{\frac{\eta}{t_0}}\,.
\end{eqnarray}
By Corollary 1.4 of \cite{P} it follows that
\begin{eqnarray}
S_q[S_u]&=&{[\,S_{{t_1}}[S_w]\,,
S_{t_0}[S_{t_0}]\,]}_\eta\,.\label{65}
\end{eqnarray}
Let $w'$ denote the conjugate of $w$\,. Note that $w'< 2< p<
{t_1}$\,. Hence, as justified above, $S_{{t_1}}[S_w]$ is
$OUMD_p$\,, while by Proposition \ref{finalsqumdp},
$S_{t_0}[S_{t_0}]\simeq S_{t_0}(\mathbb{N}\times \mathbb{N})$ is
$OUMD_p$\,.
Therefore, by (\ref{65}) it follows by interpolation that $S_q[S_u]\,$ is $OUMD_p\,.$ Furthermore, since $\max\{u', q\}\leq \max\{u', q, 2\}$ an application of (8) in Proposition \ref{propertiesumdp} shows that $S_q[S_u]$ is $OUMD_p$\,, for $\max\{u', q\}\leq p< \infty\,.$\\
To summarize, so far we have showed that $S_q[S_u]$ is $OUMD_p$\,,
for $\min\{u', q\}\leq p< \infty\,.$ It remains to analyze the
case when $1< p< \min\{u', q\}\,.$ This implies that $\min\{u,
q'\}< p'< \infty\,,$ where $q'$ is the conjugate exponent of
$q\,.$ By what we proved above it follows that $S_{q'}[S_{u'}]$ is
$OUMD_{p'}$\,. By Corollary 1.8 of \cite{P},
$(S_{q'}[S_{u'}])^*=S_q[S_u]$\,. Hence, by (5) in Proposition
\ref{propertiesumdp} it follows that $S_q[S_u]$ is $OUMD_p\,.$
This completes the proof.\qed

\begin{cor}\label{mixedcolrowha}
Let $1< p, q, u< \infty\,.$ Then the spaces $C_u\otimes^{h} R_q\,,
C_u\otimes^{h} C_q\,, R_u\otimes^{h} C_q$ and $R_u\otimes ^{h}
R_q$ are all $OUMD_p$\,.
\end{cor}
\vspace*{0.3cm} \noindent {\em Proof.} By the injectivity and
associativity properties of the Haagerup tensor product, we obtain
\begin{eqnarray*}
C_u\otimes^{h} R_q\subseteq C_u\otimes^{h}C_q\otimes^{h}
R_q\otimes^{h} R_u = S_u[S_q]\,.
\end{eqnarray*}
By (2) of Proposition \ref{propertiesumdp} and Proposition
\ref{sqsuumdp}\,, it follows that $C_u\otimes^{h} R_q$ is
$OUMD_p$\,. A similar argument applies for the spaces
$C_u\otimes^{h} C_q\,, R_u\otimes^{h} C_q\,$ and $R_u\otimes ^{h}
R_q$\,.\qed
Recall that a $C^*$-algebra $A$ has the {\em weak expectation
property} ($WEP$) of Lance \cite{La1}, if for the universal
representation $A\subset A^{**}\subset \B(H)$ there exists a
contraction $P: \B(H)\rightarrow {A^{**}}$ such that
$P\!\restriction\! A=\text{Id}_{A}\,.$ A $C^*$-algebra $B$ is said
to be $QWEP$ if it is a {\em quotient} of a $WEP$ $C^*$-algebra;
more precisely, there exists a $C^*$-algebra $A$ with the $WEP$
and a closed two-sided ideal $I$ such that $B=A/I\,.$ It is a long
standing problem whether every $C^*$-algebra is $QWEP$ (see
Kirchberg \cite{Ki} for many equivalent formulations). Note that
an injective von Neumann algebra has the $WEP$, and, therefore, it
is $QWEP$. Also, it was proved by Wassermann \cite{Wa} that for
$n\geq 2\,, VN(\mathbb{F}_n)$ is $QWEP\,,$ where $\mathbb{F}_n$ is
the free goup on $n$ generators.
% The following result is a generalization of Corollary \ref{hypnonhyp} to the case of an arbitrary $QWEP$ von Neumann algebra.

\begin{prop}\label{qwep}
Let $(\M, \tau)$ be a QWEP von Neumann algebra equipped with an
n.f. tracial state $\tau\,.$ Then $L_q(\M, \tau)$ is $OUMD_p$\,,
for $1< p, q< \infty\,.$
\end{prop}
\vspace*{0.3cm} \noindent {\em Proof}. By results of Junge
\cite{Ju}, it follows that $L_q(\M, \tau)$ is completely
contractively complemented in an ultrapower $(S_q)_{\U}\,.$
Therefore, combining Theorem \ref{finalsqumdp} with (6) and (2) of
Proposition \ref{propertiesumdp}, we obtain the conclusion.\qed

\re As a corollary, it follows that if $G$ is an amenable group
(in which case $VN(G)$ is an injective von Neumann algebra), or
$G=\mathbb{F}_n (n\geq 2)$\,, then $L_q(VN(G), \tau)$ is
$OUMD_p\,,$ for $1< p, q< \infty\,.$ \mar

\begin{cor}\label{lqluumdp}
Let $1< p, q, u< \infty\,.$ If $(\M, \tau)$ is hyperfinite von
Neumann algebra and $(\N, \phi)$ is a QWEP von Neumann algebra
equipped with n.f. tracial states $\tau$ and $\phi$\,,
respectively, then $L_q(\M; L_u(\N, \phi))$ is $OUMD_p$\,.
\end{cor}
\vspace*{0.3cm} \noindent {\em Proof.} Let $m\geq 1\,.$ Using
Junge's results from \cite{Ju}, together with the injectivity and
associativity properties of the Haagerup tensor product, we obtain
a completely contractive inclusion
\begin{eqnarray}
S_q^m[L_u(\N, \phi)]\hookrightarrow
S_q^m[(S_u)_{\U}]\,.\label{9098}
\end{eqnarray}
By Lemma 5.4 of \cite{P}, we have the complete isometry
\begin{eqnarray}
{S_q^m[(S_u)_{\U}]}&=&(S_q^m[S_u])_{\U}\,.\label{9097}
\end{eqnarray}
Therefore, by Proposition \ref{sqsuumdp}, together with (6), (1)
and (2) of Proposition \ref{propertiesumdp} it follows that
$S_q^m[L_u(\N, \phi)]$ is $OUMD_p$ and, moreover,
\begin{eqnarray*}
c_p(S_q^m[L_u(\N, \phi)])&\leq &c_p(S_q[S_u])\,.
\end{eqnarray*}
Since the spaces $S_q^m[L_u(\N, \phi)]$ are dense in $L_q(\M;
L_u(\N, \phi))$ and $OUMD_p$ is a local property, the conclusion
follows.\qed

We end this section with a discussion about the operator space
$UMD$ property for the noncommutative Lorentz spaces $L_{q, s}(\M,
\tau)$ associated to a semifinite von Neumann algebra\,. We
briefly recall some definitions. Let $x$ be a $\tau$-measurable
operator affiliated with $(\M, \tau)$\,. Following Fack and Kosaki
\cite{FK}, the {\em t-th singular number of $x$} is
\begin{eqnarray*}
\mu_t(x)&=&\inf\{\|xe\| : e \; \mbox{is a projection in } \, \M
\,,\; \tau(1-e)\leq t\,\}\,.
\end{eqnarray*}
%Recall that if $\M=L_\infty((0, \infty), m)$\,, where $m$ is the
%Lebesgue measure on $(0, \infty)$\,, and $x\in \M$ then the map
%$t\mapsto \mu_t(x)$ is exactly the {\em decreasing rearrangement}
%of the function $|x|$ (see Fack \cite{F}).
Following the general scheme of {\em symmetric operator spaces}
associated to $(\M, \tau)$ and a rearrangement invariant  Banach
function space developed in  \cite{DDP1} and \cite {DDP2}, the
noncommutative Lorentz spaces are defined as
\begin{eqnarray*}
L_{q, s}(\M, \tau)=\{x\in {L_0(\M)} :  \mu(x)\in {L_{q, s}((0,
\infty), m)}\}\,, \quad \|x\|_{L_{q, s}(\M, \tau)}=\|\mu(x)\|_{p,
q}\,.
\end{eqnarray*}
We refer the reader to Randrianantoanina \cite{Ra2} and Xu
\cite{X2, X} for more details on the noncommutative Lorentz
spaces.

\begin{prop}
Let $(\M, \tau)$ be a hyperfinite von Neumann algebra, equipped
with an n.s.f. tracial state $\tau\,.$ Then $L_{q, s}(\M, \tau)$
is $OUMD_p$\,, for $1< q, s, p< \infty$\,.
\end{prop}
 \vspace*{0.3cm}\noindent {\em Proof.} There exists $0< \theta< 1$
and $1< q_1, q_2< \infty$\,, such that
$\,\frac1q\,=\,{\frac{1-\theta}{q_0}}+{\frac{\theta}{q_2}}\,.$ The
formula describing the K-functional for the couple $(\M, L_1(\M,
\tau))$ (\cite{PX1}, Corollary 2.3) together with the reiteration
result proved in \cite{X}, Theorem 5.4(ii) yield the complete
isometry
\begin{eqnarray}
L_{q, s}(\M, \tau)&=&{[L_{q_1}(\M, \tau), L_{q_2}(\M,
\tau)]}_{\theta, s}\,.\label{53}
\end{eqnarray}
We refer to Xu's paper \cite{X} for details on the operator space
structure of the real interpolation space ${[L_{q_1}(\M, \tau),
L_{q_2}(\M, \tau)]}_{\theta, s}\,.$
%Since $\M$ is an injective von Neumann algebra, it follows by Proposition \ref{qwep} that both $\,L_{q_1}(\M, \tau)$ and %$\,L_%%{q_2}(\M, \tau)$ are UMD$_p\,.$
In particular, it is proved in \cite{X} that the following
completely isometric embedding holds
\begin{eqnarray}
{[L_{q_1}(\M, \tau), L_{q_2}(\M, \tau)]}_{\theta,
s}&\hookrightarrow &l_s(\,\{L_{q_1}(\M, \tau)+_s
{2^{-k}}L_{q_2}(\M, \tau)\}_{k}\,;
\,2^{-{k\theta}}\,)\,,\label{54}
\end{eqnarray}
where $L_{q_1}(\M, \tau)+_s {2^{-k}}L_{q_2}(\M, \tau)$ is defined
as a quotient of $\,L_{q_1}(\M, \tau){\oplus_s}
\,{2^{-k}}L_{q_2}(\M, \tau)\,,$ for $k\geq 1$\,. Then, by (2.6)'
of \cite{P}, the space $\,l_s(\,\{L_{q_1}(\M, \tau)+_s
{2^{-k}}L_{q_2}(\M, \tau)\}_{k}\,; \,2^{-{k\theta}}\,)\,$ is a
quotient of the operator space $\,l_s(\,\{L_{q_1}(\M,
\tau){\oplus_s} {2^{-k}}L_{q_2}(\M, \tau)\}_{k}\,;
\,2^{-{k\theta}}\,)\,.$ Therefore, by (2) of Proposition
\ref{propertiesumdp} it suffices to show that
$l_s(\,\{L_{q_1}(\M){\oplus_s} {2^{-k}}L_{q_2}(\M)\}_{k}\,;
\,2^{-{k\theta}}\,)$ is $OUMD_p$\,. Note that we have the complete
isometry
\begin{equation}\label{rel76560}
\,l_s(\,\{L_{q_1}(\M){\oplus_s} {2^{-k}}L_{q_2}(\M)\}_{k}\,;
\,2^{-{k\theta}}\,)=l_s(L_{q_1}(\M))\oplus_s
l_s(\{{2^{-k}}{L_{q_2}(\M)}\}_{k}\,; \,2^{-{k\theta}})\,.
\end{equation}
%Therefore, by (8) of Proposition \ref{propertiesumdp} it suffices
%to prove that both spaces $l_s(L_{q_1}(\M, \tau))$ and
%$l_s({2^{-k}}L_{q_2}(\M, \tau)\}_{k\in \mathbb{Z}}\,;
%\,2^{-{k\theta}})$ are $OUMD_p$\,.
The space $l_s(L_{q_1}(\M, \tau))$ is completely isometric to the
subspace of $S_s[L_{q_1}(\M, \tau)]$ formed by all the diagonal
matrices (see \cite{P}); moreover, the usual projection onto this
subspace is a complete contraction. Since $\M$ is hyperfinite, it
follows by Corollary \ref{lqluumdp} that $S_s[L_{q_1}(\M,
\tau)]$\, is $OUMD_p$\,. Thus, by (1) of Proposition
\ref{propertiesumdp}, it follows that $l_s(\,L_{q_1}(\M, \tau)\,)$
is $OUMD_p\,.$ Similar arguments show that
$l_s(\{{2^{-k}}{L_{q_2}}(\M, \tau)\}_{k}\,; \,2^{-{k\theta}})$ is
$OUMD_p$\,, as well. The conclusion follows from (7) of
Proposition \ref{propertiesumdp}.\qed

\vspace*{0.3cm} \noindent {\bf Question}:
%Recall that $C$ denotes the column Hilbert space of $\,l_2$ (see Section 1.2).
It was a question of Zhong-Jin Ruan whether the column Hilbert
space $C$ is $OUMD_p$\, for some (all) $1< p< \infty\,.$ By Lemma
\ref{connectumd}, if $C$ is $OUMD_p$ for some $1< p< \infty\,,$
then $S_p[C]$ is $UMD$ as a Banach space. Using the
characterization of superreflexivity in terms of ultraproducts
(see Heinrich \cite{Hei}), together with results of Junge and
Sherman from \cite{JS}, we were able to prove that $S_p[C]$ is a
superreflexive Banach space. During discussions initiated by
Gilles Pisier, Timur Oikhberg found a surprisingly short proof of
this result, which we present below. Recall that by
(\ref{cpandcpm}) we have $C_2=[C, R]_\frac12=[R,
C]_\frac12=R_2$\,. By (\ref{spehaage}) and Kouba's theorem it
follows that
\begin{equation}\label{s2c}
S_2[C]=[C, R]_\frac12\otimes^{\text{h}} C\otimes^{\text{h}} [C,
R]_\frac12=[C\otimes^{\text{h}} C\otimes^{\text{h}} C,
R\otimes^{\text{h}} C\otimes^{\text{h}} R]_\frac12\,.
\end{equation} Note that as a Banach space, $C\otimes^{\text{h}}
C\otimes^{\text{h}} C$ is isometric to a Hilbert space, hence it
is superreflexive. Pisier (see \cite{Pi8}) proved that, given a
compatible couple of Banach spaces $B_0, B_1$\,, if one of them is
superreflexive, then for $0< \theta< 1$ the interpolation space
$[B_0, B_1]_\theta$ is superreflexive, as well. Therefore, by
(\ref{s2c}) we conclude that $S_2[C]$ is a superreflexive Banach
space.
%Furthermore, the space $S_\infty[C]= C\otimes^{\text{h}}
%C\otimes^{\text{h}} R$ is superreflexive, since it is isometric
%(as a Banach space) to $K(l_2)$\,, the compact operators on
%$l_2$\,.
For $1< p< \infty$\,, by (\ref{spebyinterpo}) we have the
(complete) isometry
\[ S_p[C]=[S_\infty[C], S_2[C]]_\frac1{2p}\,, \]
A further application of Pisier's result mentioned above yields
the conclusion.

\thanks{}

\end{document}